\documentclass[a4paper,reqno,10pt]{amsart}

\usepackage[backend=bibtex,style=numeric,sorting=nyt,giveninits=true,autocite=inline]{biblatex}

\let\RNbiblatex\Rn		

\AtBeginBibliography{\let\Rn\RNbiblatex\let\em\itshape}	

\DeclareFieldFormat{labelnumberwidth}{#1.}

\DeclareFieldFormat
  [article,inbook,incollection,inproceedings,patent,thesis,unpublished]
  {title}{\mkbibitalic{#1}}

\renewbibmacro*{in:}{}

\DeclareFieldFormat{pages}{#1}

\DeclareSourcemap{
  \maps[datatype=bibtex]{
    \map{
      \step[fieldsource=mrnumber, final]
      \step[fieldset=addendum, origfieldval]
    }
  }
}
\DeclareFieldFormat{addendum}{MR~#1}

\DeclareFieldFormat{postnote}{#1}			
\DeclareFieldFormat{multipostnote}{#1}		

\DeclareSourcemap{
  \maps[datatype=bibtex]{
    \map{
      \step[fieldset=doi, null]		
      \step[fieldset=url, null]		
      \step[fieldset=urlyear, null]	
      \step[fieldset=issn, null]		
      \step[fieldset=isbn, null]		
    }
  }
}

\usepackage{geometry}				

\usepackage{comment} 

\usepackage{hyphenat}

\usepackage{manfnt} 
\reversemarginpar 

\usepackage{relsize}

\usepackage[svgnames]{xcolor}

\usepackage{url}

\usepackage{footmisc}

\usepackage{amsmath}
\usepackage{amsthm}
\usepackage{amssymb}

\allowdisplaybreaks 

\usepackage{amscd}

\usepackage{accents}

\usepackage{bbm}

\usepackage[scaled=.98]{BOONDOX-uprscr}

\usepackage[OMLmathsfit]{isomath}

\usepackage{mathabx}

\makeatletter
\DeclareSymbolFont{cmlargesymbols}{OMX}{cmex}{m}{n}
\DeclareMathSymbol{\braceld}{\mathord}{cmlargesymbols}{"7A}
\DeclareMathSymbol{\bracerd}{\mathord}{cmlargesymbols}{"7B}
\DeclareMathSymbol{\bracelu}{\mathord}{cmlargesymbols}{"7C}
\DeclareMathSymbol{\braceru}{\mathord}{cmlargesymbols}{"7D}
\DeclareRobustCommand\downbracefill{$\m@th \setbox\z@\hbox{$\braceld$}%
  \braceld\leaders\vrule \@height\ht\z@ \@depth\z@\hfill\braceru
  \bracelu\leaders\vrule \@height\ht\z@ \@depth\z@\hfill\bracerd$}
\DeclareRobustCommand\upbracefill{$\m@th \setbox\z@\hbox{$\braceld$}%
  \bracelu\leaders\vrule \@height\ht\z@ \@depth\z@\hfill\bracerd
  \braceld\leaders\vrule \@height\ht\z@ \@depth\z@\hfill\braceru$}
\DeclareRobustCommand\overbrace[1]%
     {\mathop{\vbox{\m@th\ialign{##\crcr\noalign{\kern3\p@}%
      \downbracefill\crcr\noalign{\kern3\p@\nointerlineskip}%
      $\hfil\displaystyle{#1}\hfil$\crcr}}}\limits}
\DeclareRobustCommand\underbrace[1]{\mathop{\vtop{\m@th\ialign{##\crcr
   $\hfil\displaystyle{#1}\hfil$\crcr
   \noalign{\kern3\p@\nointerlineskip}%
   \upbracefill\crcr\noalign{\kern3\p@}}}}\limits}
\makeatother

\usepackage{leftindex}

\usepackage{enumitem}
\setlist{noitemsep}

\usepackage{stackengine}
\usepackage{scalerel}

\usepackage{hyperref}
\hypersetup{colorlinks=true,linkcolor=blue,urlcolor=blue,citecolor=blue}

\newtheoremstyle{exercise}
  {3pt} 
  {3pt} 
  {\small\rmfamily} 
  {
} 
  {\rmfamily\scshape} 
  {.} 
  {.5em} 
  {} 

\newtheoremstyle{newplain}
  {5pt}
  {5pt}
  {\itshape}
  {}
  {\rmfamily\scshape}
  { ---}
  {.5em}
  {}

\newtheoremstyle{newremark}
  {5pt}
  {5pt}
  {\rmfamily}
  {}
  {\rmfamily\scshape}
  { ---}
  {.5em}
  {}

\theoremstyle{newplain}
\newtheorem*{Theorem*}{Theorem} 

\theoremstyle{newplain}
\newtheorem{Theorem}{Theorem}
\newtheorem{Lemma}[Theorem]{Lemma}
\newtheorem{Corollary}[Theorem]{Corollary}
\newtheorem{Proposition}[Theorem]{Proposition}
\newtheorem{Conjecture}[Theorem]{Conjecture}
\newtheorem{Definition}[Theorem]{Definition}

\newtheorem{IntroTheorem}{Theorem}

\theoremstyle{newremark}
\newtheorem{Empty}[Theorem]{}
\newtheorem{Remark}[Theorem]{Remark}

\newtheorem{Claim}[Theorem]{Claim}

\theoremstyle{exercise}

\numberwithin{Theorem}{section}

\newcommand{\N}{\mathbb{N}}

\newcommand{\R}{\mathbb{R}}

\newcommand{\Rm}{\R^m}
\renewcommand{\Rn}{\R^n}	

\newcommand{\bbL}{\mathbb{L}}

\newcommand{\calC}{\mathscr{C}}

\newcommand{\calD}{\mathscr{D}}

\newcommand{\calF}{\mathscr{F}}
\newcommand{\calG}{\mathscr{G}}
\newcommand{\calH}{\mathscr{H}}

\newcommand{\calS}{\mathscr{S}}
\newcommand{\calT}{\mathscr{T}}

\newcommand{\ssfL}{\mathsfit{L}}

\newcommand{\balpha}{\boldsymbol{\alpha}}

\newcommand{\bdelta}{\boldsymbol{\delta}}

\newcommand{\bkappa}{\boldsymbol{\kappa}}

\newcommand{\bmu}{\boldsymbol{\mu}}

\newcommand{\btau}{\boldsymbol{\tau}}

\newcommand{\bD}{\mathbf{D}}
\newcommand{\bE}{\mathbf{E}}
\newcommand{\bF}{\mathbf{F}}
\newcommand{\bG}{\mathbf{G}}

\newcommand{\bM}{\mathbf{M}}
\newcommand{\bN}{\mathbf{N}}

\newcommand{\bc}{\pmb{c}}

\newcommand{\bp}{\pmb{p}}
\newcommand{\bq}{\pmb{q}}

\DeclareMathOperator{\rmbdry}{\mathrm{bdry}}          

\DeclareMathOperator{\rmdiam}{\mathrm{diam}}          

\DeclareMathOperator{\rmdiv}{\mathrm{div}}            
\newcommand{\err}{\mathrm{err}}

\DeclareMathOperator{\rmev}{\mathrm{ev}}              
\DeclareMathOperator{\rmHom}{\mathrm{Hom}}            
\newcommand{\rmid}{\mathrm{id}}                       
\DeclareMathOperator{\rmint}{\mathrm{int}}            
\DeclareMathOperator{\rmLip}{\mathrm{Lip}}            

\DeclareMathOperator{\rmloc}{\mathrm{loc}}

\DeclareMathOperator{\rmosc}{\mathrm{osc}}            

\DeclareMathOperator{\rmspt}{\mathrm{spt}}            
\DeclareMathOperator{\rmSym}{\mathrm{Sym}}			
\DeclareMathOperator{\rmtrace}{\mathrm{trace}}        

\newcommand{\rmI}{\mathrm{I}}
\newcommand{\rmII}{\mathrm{II}}

\newcommand{\ip}{\,\,\begin{picture}(-1,1)(-1,-2.5)\circle*{2}\end{picture}\;\,\,}

\newcommand{\ind}{\mathbbm{1}}

\newcommand{\lseg}{\boldsymbol{[}\!\boldsymbol{[}}
\newcommand{\rseg}{\boldsymbol{]}\!\boldsymbol{]}}

\newcommand{\hel} {
\hskip2.5pt{\vrule height7pt width.5pt depth0pt}
\hskip-.2pt\vbox{\hrule height.5pt width7pt depth0pt}
\, }

\newcommand{\lno}{\pmb{\thickvert}}
\newcommand{\rno}{\mkern+1mu\pmb{\thickvert}}
\def\XXint#1#2#3{{%
\setbox0=\hbox{$#1{#2#3}{\int}$}
\vcenter{\hbox{$#2#3$}}\kern-.5\wd0}}

\usepackage{xspace}

\newcommand{\cqfd} {
\renewcommand{\qedsymbol}{$\blacksquare$}
\qed
\renewcommand{\qedsymbol}{$\blacksquare$} }
\renewcommand{\qedsymbol}{$\blacksquare$}

\newcommand{\veps}{\varepsilon}

\newcommand{\vphi}{\varphi}

\newcommand{\la}{\langle}
\newcommand{\ra}{\rangle}

\newcommand{\ie}{{\it i.e.}\ }
\newcommand{\eg}{{\it e.g.}\ }

\renewcommand{\em}{\bf}

\renewcommand{\leq}{\leqslant}
\renewcommand{\geq}{\geqslant}

\renewcommand{\subset}{\subseteq}

\makeatletter
\ifcsname phantomsection\endcsname
    \newcommand*{\qrr@gobblenexttocentry}[5]{}
\else
    \newcommand*{\qrr@gobblenexttocentry}[4]{}
\fi
\newcommand*{\addsubsection}{%
    \addtocontents{toc}{\protect\qrr@gobblenexttocentry}%
    \subsection}
\makeatother

\begin{document}



\title[]{The divergence of continuously differentiable vector fields}

\author[Th. De Pauw]{Thierry De Pauw}

\address{Institute for Theoretical Sciences / School of Science, Westlake University\\
No. 600, Dunyu Road, Xihu District, Hangzhou, Zhejiang, 310030, China}

\email{thierry.depauw@westlake.edu.cn}




\begin{abstract}
For $m\geq 2$, we give a complete characterisation of the continuous functions on $\Rm$ that are divergences of continuously differentiable vector fields.
The characterisation is formulated in terms of distributions and, equivalently, in terms of the continuous dual of a suitably topologised subspace of the Arens--Eells space.
The two implications yield contrasting applications.
First, there exists a compactly supported continuous function which is the divergence of a bounded Lipschitz vector field but not of any continuously differentiable vector field.
Second, every continuous function whose distributional partial derivatives are locally $m$-integrable (the critical exponent for Sobolev embedding) is the divergence of a continuously differentiable vector field.
\end{abstract}

\maketitle



\tableofcontents


\section{Introduction}

If $v : \Rm \to \Rm$ is continuously differentiable, then its divergence $f = \rmdiv v$ is continuous.
Conversely, given a continuous function $f : \Rm \to \R$, does there exist a continuously differentiable vector field $v$ such that $\rmdiv v = f$?
\par
If $m=1$, the answer is positive.
Indeed, Riemann integration defines $I(f)(x)=\int_0^x f(t)\,dt$.
Thus $I:C(\R)\to C^1(\R)$ is an explicit, linear and continuous right inverse of differentiation, the two spaces being equipped with their usual Fr\'echet topologies.
For $m\geq 2$, integration in a single coordinate direction still provides a solution of the divergence equation, but it need not produce a $C^1$ vector field.
\par
The classical approach is to solve the Poisson equation $\triangle g=f$ and put $v=\nabla g$.
For instance, when $f$ has compact support, one may use the Newtonian potential $g=\Gamma*f$, yielding
\begin{equation*}
v=\nabla(\Gamma*f) \quad\text{and}\quad \rmdiv v=f.
\end{equation*}
The Calder\'on--Zygmund theory gives $L_p$ estimates for $Dv$ when $f\in L_p$ and $1<p<\infty$ (see \eg \cite{STEIN.70}).
H\"older regularity of $f$ yields $v\in C^{1,\alpha}$ when $f\in C^\alpha$, $0<\alpha<1$ (see \eg \cite{GILBARG.TRUDINGER}).
These are quantitative regularity statements.
A Dini modulus of continuity of $f$ also suffices for continuity of $Dv$, whereas mere continuity does not (see \cite[Problems 4.2 and 4.9(a)]{GILBARG.TRUDINGER}).
Similarly, $f\in L_\infty$ yields $Dv\in\mathrm{BMO}$, but not in general $Dv\in L_\infty$, see \eg McMullen \cite{MCM.98}.
\par
There is, however, a difference between the Laplacian and the divergence operator.
Two distributional solutions of $\triangle g=f$ differ by a harmonic function, hence by a real-analytic function.
Thus, for the Poisson equation, investigating the regularity of one solution suffices to decide the regularity of all solutions.
The kernel of divergence, on the contrary, contains irregular vector fields.
A solution of the divergence equation may, therefore, be considerably smoother than the one supplied by the Newtonian potential.
\par
Notwithstanding, Preiss \cite{PRE.97} proved that, when $m\geq 2$, there exists $f\in C_c(\Rm)$ for which $\rmdiv v=f$ has no locally Lipschitz solution.
In particular, no $C^1$ solution exists.
His proof relies on Lipschitz determining sets of measure zero and on the construction of an oscillating function.
McMullen \cite{MCM.98} obtained another obstruction for bounded measurable data by a duality argument.
For a recent treatment of these topics, see Curc\u a \cite{CUR.22}.
Dacorogna, Fusco and Tartar \cite{DAC.FUS.TAR.03} also established the failure of surjectivity for continuous data and considered the $L_1$ endpoint.
For an overview of divergence equations, including boundary conditions, see Russ \cite{RUS.13}.
\par
The arguments of McMullen and of Dacorogna, Fusco and Tartar make use of Ornstein's non-inequality \cite{ORN.62}.
Recall that the Calder\'on--Zygmund theory gives, for $1<p<\infty$, an estimate
\begin{equation*}
\|\triangle\vphi\|_{L_p}\leq \bc(m,p)\cdot\|D_0^2\vphi\|_{L_p},
\end{equation*}
$\vphi\in\calD(\Rm)$, where $D_0^2\vphi=D^2\vphi-\frac{1}{m}(\triangle\vphi)I_m$ is the trace-free Hessian.
Ornstein's theorem asserts, in particular, that no corresponding estimate holds for $p=1$.
The form of this failure used in our construction is recalled in \ref{Ornstein}.
The preceding counterexamples answer the existence question negatively, but do not characterise the continuous functions which are divergences of $C^1$ vector fields.
This is our objective.
The characterisation we offer relies on the interplay between linear partial differential equations, abstract functional analysis and geometric measure theory.
\par
To state our first result, we let $\calD_\#(\Rm)$ denote the space of test functions whose integral vanishes:
\begin{equation*}
\calD_\#(\Rm)=\calD(\Rm) \cap \left\{u \mid \int_{\Rm}u\,d\ssfL^m=0\right\}.
\end{equation*}
For $u\in\calD_\#(\Rm)$, we define
\begin{equation*}
\bG(u)=\sup\left\{\int_{\Rm}u\cdot\vphi\,d\ssfL^m \mid \vphi\in\calD(\Rm) \text{ and } \|\nabla\vphi\|_\infty\leq1\right\},
\end{equation*}
and, using the same norm for each partial derivative $D_iu$, we put
\begin{equation*}
\lno u\rno=\sum_{i=1}^m\bG(D_iu).
\end{equation*}
Both are norms on $\calD_\#(\Rm)$, see \ref{3.1}.
The origin of these two norms will be explained below.
\begin{IntroTheorem}[Characterisation of $C^1$ divergences, see \ref{6.8}]
Let $F:\calD_\#(\Rm)\to\R$ be linear.
The following are equivalent.
\begin{enumerate}
\item[(A)] There exists $v\in C^1(\Rm;\Rm)$ such that
\begin{equation*}
\la u,F\ra=\int_{\Rm}u\cdot\rmdiv v\,d\ssfL^m,
\end{equation*}
for every $u\in\calD_\#(\Rm)$.
\item[(B)] For every $\veps>0$ and every $k\in\N_+$ there exists $\theta>0$ such that
\begin{equation*}
|\la u,F\ra|\leq\theta\cdot\bG(u)+\veps\cdot\lno u\rno,
\end{equation*}
for every $u\in\calD_\#(\Rm)$ supported in $B(0,k)$.
\end{enumerate}
\end{IntroTheorem}
Theorem \ref{6.8} also gives an equivalent formulation in terms of continuous linear extensions to a suitably topologised space of distributions.
For $F$ obtained by integration against a continuous function $f$, conclusion (A) gives $\rmdiv v=f$ up to an additive constant, which is easily absorbed by a linear vector field.
An early indication that the stronger norm is adapted to the problem is the estimate
\begin{equation*}
\lno\triangle\vphi\rno\leq\frac{m}{m-1}\cdot\|D_0^2\vphi\|_{L_1},
\end{equation*}
for $\vphi\in\calD(\Rm)$, established in Theorem \ref{4.2F}.

\par
To understand the choice of our two norms, we return to the definition of distributional divergence.
For a locally integrable vector field $v$ and $\vphi\in\calD(\Rm)$, it reads
\begin{equation*}
\la\vphi,\rmdiv v\ra=-\sum_{i=1}^m\int_{\Rm}(D_i\vphi)\cdot v_i\,d\ssfL^m.
\end{equation*}
The components $v_i$ are paired with the derivatives $D_i\vphi$, rather than with $\vphi$ itself.
This suggests considering distributions $T$ whose derivatives $D_iT$ can be paired with the components of $v$.
Our discussion is local: test objects have support in a prescribed compact set, and the regularity of $v$ is considered in a neighbourhood of that set.
\par
For continuous vector fields, a natural requirement is that each $D_iT$ be a signed Borel measure.
For integrable vector fields, we instead require that each $D_iT$ be represented by an essentially bounded function.
These two conditions lead to functions of bounded variation and to Lipschitz functions, respectively; see \cite[Section 10]{DEP.26c} and \cite{DEP.26a}.
When $v$ is of class $C^1$, its components are locally Lipschitz.
This suggests requiring $\bG(D_iT)<\infty$ for each $i$.
For a compactly supported distribution $T$, we denote by $\chi(T)$ its augmentation, namely its value on any test function equal to $1$ on a neighbourhood of its support.
For every $k\in\N_+$, we define
\begin{equation*}
\begin{aligned}
X_k&=\calD_0(\Rm) \cap \left\{T \mid \rmspt T\subset B(0,k),\ \chi(T)=0, \text{ and } \lno T\rno<\infty\right\},\\
X&=\bigcup_{k\in\N_+}X_k,
\end{aligned}
\end{equation*}
where $\calD_0(\Rm)$ is the space of distributions in $\Rm$.
These spaces $X$ and $X_k$ are studied in \ref{4.2}.
The nature of these test objects is less evident, however.
In particular, they need not be functions, as we shall discuss below, see \ref{5.8}.

\par
The two norms in Theorem A arise naturally from an approximation argument.
Fix a compact set $K\subset\Rm$ and a vector field $v$ defined on a neighbourhood of $K$.
Choose a smooth vector field $w$ approximating $v$ there, in the regularity appropriate to the class under consideration for $v$.
For $u\in\calD_\#(\Rm)$ supported in $K$, integration by parts yields
\begin{equation*}
\la u,\rmdiv v\ra=\int_{\Rm}u\cdot\rmdiv w\,d\ssfL^m
-\sum_{i=1}^m\int_{\Rm}(D_iu)\cdot(v_i-w_i)\,d\ssfL^m.
\end{equation*}
The first term pairs $u$ with a smooth function.
The second pairs $D_iu$ with a function that can be made small in the regularity prescribed for $v_i$.
\par
For continuous $v$, uniform approximation on $K$ gives two norms: $\|u\|_{L_1}$ and $\sum_i\|D_iu\|_{\mathrm M}$, where $\|\cdot\|_{\mathrm M}$ denotes the total variation norm of signed Borel measures.
For integrable $v$, approximation in $L_1$ gives instead $\|u\|_{L_\infty}$ and $\sum_i\|D_iu\|_{L_\infty}$.
For $v\in C^1$, the Lipschitz constants of $v_i-w_i$ can be made arbitrarily small on a neighbourhood of $K$.
After a suitable cutoff, the second term is, therefore, bounded by $\veps\cdot\lno u\rno$.
Since $\rmdiv w$ is smooth, the first term is bounded by $\theta\cdot\bG(u)$, for some $\theta>0$.
We obtain, for each $k\in\N_+$ and $\veps>0$, an estimate
\begin{equation*}
|\la u,\rmdiv v\ra|\leq\theta\cdot\bG(u)+\veps\cdot\lno u\rno,
\end{equation*}
for all $u\in\calD_\#(\Rm)$ supported in $B(0,k)$.
This explains the two norms and the form of the necessary condition in Theorem A.
The converse requires more than mollification, as we shall see below.

\par
The conditions imposed upon the derivatives of the test objects also tell us something about the test objects themselves.
If the distributional derivatives $D_iT$ are signed Borel measures, then $T$ is represented, locally, by a function of bounded variation.
If instead $D_iT\in L_\infty$ for all $i$, then $T$ is represented, locally, by a Lipschitz function.
In either case, the test objects are functions, up to the usual choice of an additive constant.
\par
The corresponding conclusion for $X$ is less familiar.
Theorem \ref{4.5} shows that $\bG(T)<\infty$ for $T\in X$.
Moreover, the mollification estimate in \ref{4.7}(D) implies that $T$ is a flat $0$-chain, see \ref{5.3}.
Using the identification recalled in \ref{np.arens-eells}, we obtain, for every $k\in\N_+$,
\begin{equation*}
X_k\subset\calG(B(0,k))\cong\calF(B(0,k)),
\end{equation*}
where $\calF(B(0,k))$ is the Arens--Eells, or Lipschitz-free, space of $B(0,k)$, see \ref{5.4}.
One may think of this space as the completion, in the $\bG$ norm, of the finite linear combinations of differences of Dirac masses.
\par
Notice that the membership of $T$ in an Arens--Eells space does not imply the corresponding membership of its derivatives $D_iT$, see \ref{5.6}.
This is analogous to the fact that a function of bounded variation belongs to $L_1$, while its distributional derivatives need not belong to $L_1$.
More strikingly, the members of $X$ need not even be represented by signed Borel measures, as shown in \ref{5.8}.
Thus the use of distributions, rather than functions, is intrinsic to the present test space.

\par
A further analogy concerns Poincar\'e-type inequalities.
For functions of bounded variation supported in a fixed bounded set, the $L_1$ norm is controlled by the total variation of their distributional gradient.
For Lipschitz functions with the same support restriction, the $L_\infty$ norm is controlled by the $L_\infty$ norm of their gradient.
The corresponding result for $X$ is Theorem \ref{4.5}: for each $k\in\N_+$ there is a constant $\bc(m,k)>0$ such that
\begin{equation*}
\bG(T)\leq\bc(m,k)\cdot\lno T\rno
\end{equation*}
whenever $T\in X_k$.
In each case, the weaker norm of the test object is thus controlled by the norm of its derivatives, under a restriction on its support.
\par
These inequalities explain part of the relationship between the two norms, but they do not provide the compactness needed for our existence argument.
The next question is whether, under a support restriction, families bounded in the stronger norm are relatively compact with respect to the weaker norm.
\par
The answer is affirmative in all three cases.
For functions of bounded variation, a bound on the total variation of the gradient and a uniform support restriction yield compactness in $L_1$.
This is a consequence of the Fr\'echet--Kolmogorov--Riesz--Weil criterion: the gradient bound controls translations in $L_1$, while the support restriction ensures tightness.
For Lipschitz functions, the corresponding conclusion in the uniform norm follows from Ascoli's theorem, since the bound on the gradient gives equicontinuity.
\par
For our test distributions, the analogous result is Theorem \ref{6.2}.
We put
\begin{equation*}
C_k=X_k\cap\left\{T \mid \lno T\rno\leq k\right\},
\end{equation*}
$k\in\N_+$.
The theorem asserts that each $C_k$ is compact with respect to $\bG$.
The main step of its proof is the mollification estimate
\begin{equation*}
\bG(T-T*\Phi_\delta)\leq\delta\cdot\lno T\rno,
\end{equation*}
established in \ref{4.7}(D).
For fixed $\delta>0$, the mollifications of members of $C_k$ have uniformly bounded supports and form a bounded, equi-Lipschitz family of smooth functions, by \ref{4.7}(A)--(B).
Ascoli's theorem then gives total boundedness for $\bG$, and completeness is established separately.
\par
Here, too, one may distinguish between a general compactness criterion and a sufficient condition involving derivatives.
Our earlier work \cite[theorem 8.3]{DEP.26b} characterises $\bG$-compact families in Lipschitz-free spaces by uniform approximation with objects of finite mass.
The bound on $\lno T\rno$ ensures that this criterion applies to the sets $C_k$, as explained in Remark \ref{6.3}.
Thus Theorem \ref{6.2} is a compactness result tailored to the divergence problem.
\par
The compactness of the $C_k$'s leads us to the localised locally convex topologies developed in \cite{DEP.26c}.
Let $\calG$ be the topology induced by $\bG$ on $X$, and let $\calG_{\lno\cdot\rno}$ be its localisation by the family $\calC=\{C_k \mid k\in\N_+\}$, see \ref{6.1}.
This topology induces $\calG$ on each $C_k$.
Moreover, a linear form $F:X\to\R$ is $\calG_{\lno\cdot\rno}$-continuous if and only if, for every $k\in\N_+$ and $\veps>0$, there is $\theta>0$ such that
\begin{equation*}
|\la T,F\ra|\leq\theta\cdot\bG(T)+\veps\cdot\lno T\rno
\end{equation*}
for all $T\in X_k$, compare \ref{np.loc-top}.
This recovers the form of the inequality in Theorem A.
By the compactness of the $C_k$'s, the space $X[\calG_{\lno\cdot\rno}]$ is semireflexive, according to \cite[theorem 7.4]{DEP.26c}.
In other words, its canonical map into its strong bidual is surjective.
\par
To turn this observation into an existence result, we consider the space of $C^1$ vector fields vanishing at the origin, see \ref{6.4}.
The topology of local $C^1$ convergence makes this a Fr\'echet space.
The distributional pairing extends, by approximation, to a linear map from this space to the continuous dual of $X[\calG_{\lno\cdot\rno}]$, see \ref{6.6}.
In \ref{6.7}, we verify the hypotheses of the abstract existence theorem \cite[theorem 8.1]{DEP.26c}.
Compactness alone does not give surjectivity.
The proof in \cite{DEP.26c} also uses the Hahn--Banach and closed range theorems, together with the Krein--\v{S}mulian theorem.
This yields the sufficient implication of Theorem \ref{6.8} and completes the characterisation.
\par
Ideas related to localised topologies were discussed earlier and in a different context by the author in \cite{DEP.97}.
Later, motivated by the work of Bourgain and Brezis \cite{BOU.BRE.03}, we subsequently developed this approach in joint work with Pfeffer \cite{DEP.PFE.06b} and with Moonens and Pfeffer \cite{DEP.MOO.PFE.08}.
The abstract theory we recently set forth in \cite{DEP.26c} places these methods in a common functional-analytic setting.
The present application shows that, once the appropriate test distributions and their compactness have been identified, this setting also applies to $C^1$ vector fields.
\par
We now present two applications of our characterisation.
The first one distinguishes $C^1$ solvability from Lipschitz solvability.
\begin{IntroTheorem}[A Lipschitz divergence without $C^1$ solution, see \ref{8.4}]
There exists $f\in C_c(\Rm)$ for which $\rmdiv v=f$ admits a bounded Lipschitz solution $v:\Rm\to\Rm$, but no solution of class $C^1$.
\end{IntroTheorem}
This conclusion is different from the obstructions identified by Preiss and McMullen recalled above.
Indeed, it concerns a continuous function which is already known to be the divergence of a Lipschitz vector field.
\par
The proof begins with a dual form of Ornstein's non-inequality, Theorem \ref{8.1}.
It leads to smooth data with small supremum norm, but whose pairing with suitable test functions remains bounded away from zero, see \ref{8.3}.
Translating and rescaling these examples, we construct a compactly supported continuous function $f$ and a bounded Lipschitz vector field with divergence $f$.
If $f$ were also the divergence of a $C^1$ vector field, the necessary estimate of Theorem A, applied to the rescaled test functions, would give a contradiction.
The argument is based on Ornstein's original non-inequality.
For related formulations and proofs, see Kirchheim--Kristensen \cite{KIR.KRI.16}, Conti--Faraco--Maggi \cite{CON.FAR.MAG.05}, and Faraco--Guerra \cite{FAR.GUE.22}.
\par
Our second application gives a positive result for continuous functions having first distributional derivatives in $L_{m,\rmloc}$.
\begin{IntroTheorem}[A sufficient condition for $C^1$ solvability, see \ref{9.4}]
Let $f\in C(\Rm)$ and assume that $D_i f\in L_{m,\rmloc}(\Rm)$ for all $i\in\{1,\ldots,m\}$.
Then there exists $v\in C^1(\Rm;\Rm)$ such that $\rmdiv v=f$.
\end{IntroTheorem}
The exponent $m$ is critical for the Sobolev embedding.
The condition on the derivatives is weaker than the Lorentz-space condition $\nabla f\in L_{m,1,\rmloc}$, which, by Stein's endpoint theorem \cite{STE.81}, ensures that $f$ has a continuous representative.
Here, continuity is assumed separately, since $W^{1,m}_{\rmloc}$ alone does not imply it.

\par
The principal estimate behind Theorem C is established in \ref{9.3}.
For $u\in\calD_\#(\Rm)$ supported in $B(0,k)$, it gives
\begin{equation*}
\left|\int_{\Rm}u\cdot f\,d\ssfL^m\right|\leq
\bc(m,k)\cdot\left(\|f\|_{\infty,B(0,k+2)}+\|\nabla f\|_{L_m,B(0,k+2)}\right)\cdot\lno u\rno.
\end{equation*}
To see where the two terms on the right come from, we first represent each $D_i u$ as the divergence of a smooth vector field with $L_1$ norm controlled by $\bG(D_i u)$, see \ref{9.2}.
These fields are assembled into a matrix field $A$.
With $g$ so that $\triangle g=f$ near the support of $A$, integration by parts expresses the pairing $\int u \cdot f \,d\ssfL^m$ as $\int A:D^2g\,d\ssfL^m$.
The trace part of $A$ is controlled by $\|f\|_\infty$, since $\triangle g=f$.
For the trace-free part, the differential constraint on $A$ allows us to apply Van Schaftingen's estimate for cocancelling operators \cite{VSJ.13}.
The remaining derivatives of $g$ are estimated in $L_m$ by Calder\'on--Zygmund theory.
\par
We next approximate $f$ by its mollifications.
In the decomposition of $\int u\cdot f\,d\ssfL^m$ into a smooth contribution and a remainder, the former is controlled by $\theta\cdot\bG(u)$.
The estimate from \ref{9.3}, together with local uniform convergence of the mollifications and $L_m$ convergence of their gradients, makes the latter smaller than $\veps\cdot\lno u\rno$.
The sufficient condition of Theorem A then applies, proving Theorem C.
%
%
%
%
\par
We return to a feature of the one-dimensional case: the dependence of a solution upon the datum.
Riemann integration and the Newtonian potential give explicit, continuous linear solution operators in the settings where they apply.
For the divergence operators associated with continuous and integrable vector fields, by contrast, no continuous linear right inverse exists, see \cite[Section 1.4]{DEP.26a}, where the results of Bouafia \cite{BOU.11} and of Naor--Schechtman \cite{NAO.SCH.07} are recalled.
Our approach to the $C^1$ case does not provide a linear right inverse, though we conjecture that it does not exist.
The abstract theory \cite[Theorem 8.1(C4)]{DEP.26c} does, however, provide continuous, in general nonlinear, selections with quantitative control, by an application of Michael's selection theorem.
In particular, Theorem \ref{9.4} gives the following estimate: for each $k\in\N_+$ such that $f$ is not constant on $B(0,k)$, there exists a solution $v$ of $\rmdiv v=f$, possibly depending on $k$, such that
\begin{equation*}
\bq_k(v)\leq\bc(m,k)\cdot\left(\|f\|_{\infty,B(0,k+3)}+\|\nabla f\|_{L_m,B(0,k+3)}\right).
\end{equation*}
Here, $\bq_k$ measures the size of the first derivatives of $v$ on $B(0,k)$, see \ref{np.omega}.
\par
The paper is organised as follows.
Sections 2--5 develop the test spaces and their estimates, including the necessary condition for $C^1$ solvability.
Section 6 establishes the characterisation of Theorem A.
Section 7 discusses Lipschitz divergences, while Sections 8 and 9 prove Theorems B and C, respectively.

\section{Notation and preliminaries}

\begin{Empty}[Euclidean notation and test functions]
\label{np.euclidean}
Throughout, $m\geq 2$.
We denote by $|\cdot|_2$ the Euclidean norm on $\Rm$ and by $a\ip b$ the Euclidean scalar product of $a,b\in\Rm$.
We let $\calD(\Rm)$ be the space of smooth compactly supported real-valued functions on $\Rm$, and put
\begin{equation*}
\calD_\#(\Rm)=\calD(\Rm) \cap \left\{u \mid \int_{\Rm}u\,d\ssfL^m=0\right\}.
\end{equation*}
We use the letter $\vphi$ for members of $\calD(\Rm)$ and the letter $u$ for members of $\calD_\#(\Rm)$ when they occur as test functions.
We write $\calD(\Rm;\Rm)$ for the space of smooth compactly supported vector fields.
For a differentiable scalar function $\vphi$, we use
\begin{equation*}
\|\nabla\vphi\|_\infty=\sup_{x\in\Rm}|\nabla\vphi(x)|_2.
\end{equation*}
The norm on matrices is introduced separately below.
\par 
For the purpose of mollifying, we fix an even, nonnegative $\Phi\in\calD(\Rm)$ supported in $B(0,1)$ such that $\int_{\Rm}\Phi\,d\ssfL^m=1$, and put $\Phi_\delta(x)=\delta^{-m} \cdot \Phi(x/\delta)$ for $\delta>0$.
\end{Empty}

\begin{Empty}[Currents, mass, boundary, and augmentation]
\label{np.currents}
We use the terminology and notation of geometric measure theory, referring to Federer \cite{GMT} and Simon \cite{SIMON} for general definitions and results.
For $k=0,\ldots,m$, we denote by $\calD^k(\Rm)$ the vector space of smooth compactly supported differential forms of degree $k$ and by $\calD_k(\Rm)$ the vector space of $k$-dimensional currents.
In particular, $\calD^0(\Rm)=\calD(\Rm)$ is the usual space of test functions, whereas $\calD_0(\Rm)$ is the space of distributions.
We write $\bF_k(\Rm)$ and $\bN_k(\Rm)$ for the spaces of $k$-dimensional flat chains and normal currents, respectively.
Mass and flat norm are denoted by $\bM$ and $\bF$.
The boundary of a current $T\in\calD_k(\Rm)$, $k\geq1$, is the current $\partial T\in\calD_{k-1}(\Rm)$ defined by
\begin{equation*}
\la\omega,\partial T\ra=\la d\omega,T\ra,
\end{equation*}
$\omega\in\calD^{k-1}(\Rm)$.
For $u\in L_1(\Rm)$, we denote by $\ssfL^m\hel u$ the associated $0$-current, defined by
\begin{equation*}
\la\vphi,\ssfL^m\hel u\ra=\int_{\Rm}\vphi\cdot u\,d\ssfL^m,
\end{equation*}
$\vphi\in\calD(\Rm)$.
Similarly, for $v\in L_1(\Rm;\Rm)$, the $1$-current $\ssfL^m\wedge v$ is defined by
\begin{equation*}
\la\omega,\ssfL^m\wedge v\ra=\int_{\Rm}\omega(x)\ip v(x)\,d\ssfL^m(x),
\end{equation*}
$\omega\in\calD^1(\Rm)$.
In particular,
\begin{equation*}
\bM(\ssfL^m\hel u)=\|u\|_{L_1}\qquad\text{and}\qquad \bM(\ssfL^m\wedge v)=\int_{\Rm}|v|_2\,d\ssfL^m.
\end{equation*}
If $v\in C_c^1(\Rm;\Rm)$, then $\partial(\ssfL^m\wedge v)=-\ssfL^m\hel\rmdiv v$.
Given a Lipschitz map $f:\Rm\to\Rm$ and a compactly supported flat chain $T$, we denote its pushforward by $f_\#T$.
We refer to Federer and Simon for the definition and properties of Lipschitz pushforwards, including their compatibility with boundary.
The pullback of a differential form is denoted by $f^\#\omega$ whenever it is defined.
For a compactly supported distribution $T\in\calD_0(\Rm)$, its {\em augmentation} $\chi(T)$ is defined by
\begin{equation*}
\chi(T)=\la\vphi,T\ra,
\end{equation*}
where $\vphi\in\calD(\Rm)$ equals $1$ on a neighbourhood of $\rmspt T$.
The value is independent of the choice of $\vphi$, and $\chi$ is linear.
In particular, $\chi(D_iT)=0$ for compactly supported $T$, and $\chi(\ssfL^m\hel u)=\int_{\Rm}u\,d\ssfL^m$ for compactly supported $u\in L_1(\Rm)$.
More generally, if a compactly supported distribution $T$ is represented by a signed Borel measure $\mu$, then $\chi(T)=\mu(\Rm)$.
\end{Empty}

\begin{Empty}[Arens--Eells spaces and flat chains]
\label{np.arens-eells}
Let $Q=[-1,1]^m$ be the unit cube.
We denote by $\rmLip_0(Q)$ the Banach space of real-valued Lipschitz functions on $Q$ vanishing at the origin, equipped with the Lipschitz constant as norm.
For every $x\in Q$, let $\delta_x\in\rmLip_0(Q)^*$ be the evaluation functional defined by
\begin{equation*}
\la f,\delta_x\ra=f(x),
\end{equation*}
$f\in\rmLip_0(Q)$.
Thus $\delta_0=0$ in $\rmLip_0(Q)^*$.
The {\em Arens--Eells space}, also called the {\em Lipschitz-free space}, $\calF(Q)$ is the closed vector subspace of $\rmLip_0(Q)^*$ generated by the evaluation functionals $\delta_x$, $x\in Q$, equipped with the induced dual norm.
In particular, for $T\in\calF(Q)$,
\begin{equation*}
\|T\|_{\calF(Q)}=\sup\left\{|\la f,T\ra| \mid f\in\rmLip_0(Q)\text{ and }\rmLip f\leq1\right\}.
\end{equation*}
The dual of $\calF(Q)$ is canonically and isometrically identified with $\rmLip_0(Q)$.
In \cite{DEP.26b}, we established an isometric isomorphism between $\calF(Q)$ and
\begin{equation*}
\calG(Q)=\bF_0(\Rm)\cap\{T \mid \rmspt T\subset Q\text{ and }\chi(T)=0\},
\end{equation*}
equipped with the norm $\bG$ defined in \ref{4.2} below.
The isomorphism sends $\delta_x\in\rmLip_0(Q)^*$ to the flat $0$-chain $\bdelta_x-\bdelta_0$, where the latter Dirac masses are understood as currents.
The same identification holds with a closed Euclidean ball $B$ containing the origin in place of $Q$, using the corresponding spaces $\calF(B)$ and $\calG(B)$, see also \cite{DEP.26a}.
In particular, the spaces $X_k$ introduced in \ref{4.2} satisfy
\begin{equation*}
X_k\subset\calG(B(0,k))\cong\calF(B(0,k)).
\end{equation*}
The inclusion is strict, as shown in \ref{5.5}.
\end{Empty}

\begin{Empty}[The balls $B(0,k)$, the $\psi_k$, and the $\bq_k$]
\label{np.omega}
For all $k \in \N_+$ we fix $\psi_k \in \calD(\Rm)$ such that $0 \leq \psi_k \leq 1$, $\psi_k = 1$ on a neighbourhood of $B(0,k)$, and $\rmspt \psi_k \subset U(0,k+1)$.
We may and shall assume that $\|\nabla \psi_k \|_\infty \leq 2$ for all $k \in \N_+$.
\par 
For $v \in C^1(\Rm;\Rm)$ and $k \in \N_+$ we define
\begin{equation*}
\bq_k(v) = \max_{i \in \{1,\ldots,m\}} \sup_{x \in B(0,k)} |\nabla v_i(x)|_2 .
\end{equation*}
Note that $\bq_k \leq \bq_{k+1}$ for all $k \in \N_+$.
\end{Empty}

\begin{Empty}[Matrices and their norms]
\label{np.matrices}
We let $M_m(\R)$ denote the vector space of real $m \times m$ matrices, $\rmSym(m)$ its subspace of symmetric matrices, and $I_m$ the identity matrix.
Given $A \in M_m(\R)$ we denote its $i^{th}$ row by $A_i$ and define
\begin{equation*}
|A| = \sum_{i=1}^m |A_i|_2,
\end{equation*}
where the norm on the right is the Euclidean norm on $\Rm$.
This is the norm we use on $M_m(\R)$, unless explicitly stated otherwise.
For $A,B \in M_m(\R)$ we put
\begin{equation*}
A:B = \sum_{i,j=1}^m A_{ij} \cdot B_{ij} = \rmtrace(A^t B).
\end{equation*}
Notice that
\begin{equation*}
|A:B| \leq |A| \cdot \max_{i} |B_i|_2.
\end{equation*}
We also record
\begin{equation*}
|\rmtrace A| \leq |A| \qquad\text{and}\qquad |I_m| = m.
\end{equation*}
The trace-free part of $A$ is
\begin{equation*}
A^\circ = A - \frac{\rmtrace A}{m} \cdot I_m.
\end{equation*}
In particular,
\begin{equation*}
\rmtrace A^\circ = 0 \qquad\text{and}\qquad |A^\circ| \leq 2 \cdot |A|.
\end{equation*}
We write
\begin{equation*}
M_{m,0}(\R) = M_m(\R) \cap \{ A \mid \rmtrace A = 0 \}.
\end{equation*}
Given $a,b \in \Rm$, their tensor product $a \otimes b \in M_m(\R)$ is defined by $(a \otimes b)_{ij} = a_i \cdot b_j$.
Thus,
\begin{equation*}
(a \otimes b) \cdot x = (b \ip x) \cdot a \qquad\text{and}\qquad A:(a \otimes b) = a \ip (A \cdot b).
\end{equation*}
If $A$ is symmetric, the sum of the Euclidean norms of its columns equals $|A|$.
When necessary, we distinguish the operator norm
\begin{equation*}
|A|_{\mathrm{op}} = \sup_{|x|_2 \leq 1} |A \cdot x|_2
\end{equation*}
from our standard matrix norm.
In particular, the Hessian bounds for radial cutoffs will be expressed using the operator norm.
\end{Empty}

\begin{Empty}[Localised topologies and an abstract existence theorem]
\label{np.loc-top}
We recall the setting and those conclusions of \cite{DEP.26c} that will be used below.
Let $X$ be a real vector space equipped with a norm $\bG$, and let $\lno \cdot \rno$ be another norm on $X$ which is lower-semicontinuous with respect to $\bG$.
Assume that $\la X_k \ra_{k \in \N_+}$ is a non-decreasing sequence of $\bG$-closed vector subspaces whose union is $X$.
Put
\begin{equation*}
C_k = X_k \cap \{ x \mid \lno x \rno \leq k \} \quad\text{and}\quad \calC = \{ C_k \mid k \in \N_+ \}.
\end{equation*}
The family $\calC$ is localising in the sense of \cite[1.2]{DEP.26c}.
We denote by $\calG$ the topology induced by $\bG$ and by $\calG_{\calC}$ its localisation by $\calC$, defined in \cite[1.3--1.4]{DEP.26c}.
The topology $\calG_{\calC}$ is characterised by the fact that it induces $\calG$ on each $C_k$ and that a linear map from $X$ into a locally convex topological vector space is $\calG_{\calC}$-continuous whenever its restrictions to the $C_k$'s are $\calG$-continuous.
According to \cite[3.1]{DEP.26c}, a sequence in $X$ converges with respect to $\calG_{\calC}$ if and only if it converges with respect to $\bG$ and is contained in some $C_k$; a subset of $X$ is $\calG_{\calC}$-bounded if and only if it is $\bG$-bounded and contained in some $C_k$.
If each $C_k$ is $\bG$-compact, then $X[\calG_{\calC}]$ is sequential and semireflexive, by \cite[4.4(C1) and 7.4]{DEP.26c}.
For clarity, recall that the dual $X[\calG_{\calC}]^*$ consists of the $\calG_{\calC}$-continuous linear forms on $X$.
Its {\em strong topology} $\beta(X^*,X)$ is the locally convex topology of uniform convergence on the $\calG_{\calC}$-bounded subsets of $X$: it is generated by the seminorms
\begin{equation*}
F \mapsto \sup_{x \in B} |\la x,F \ra|, 
\end{equation*}
corresponding to all $\calG_{\calC}$-bounded sets $B \subset X$.
The {\em strong bidual} $X^{**}$ is the continuous dual of $X^*[\beta(X^*,X)]$.
The evaluation map $\rmev : X \to X^{**}$ is defined by $\la F,\rmev(x)\ra=\la x,F\ra$.
We say that $X[\calG_{\calC}]$ is {\em semireflexive} if this map is surjective.
As the topology is Hausdorff and locally convex, it is also injective.
\par
We recall the following normed version of the abstract existence theorem \cite[Theorem 8.1]{DEP.26c}.
Assume each $C_k$ is $\bG$-compact.
Let $E$ be a Fr\'echet space whose topology is generated by a non-decreasing sequence of seminorms $\la q_k \ra_k$, and let $D:E\to X^\star$ be linear, where $X^\star$ denotes the algebraic dual.
Assume the following conditions.
\begin{enumerate}
\item[(H1)] For every $v\in E$ and $\veps>0$ there exists $\theta>0$ such that
\begin{equation*}
|\la x,D(v)\ra|\leq \theta\cdot\bG(x)+\veps\cdot\lno x\rno,
\end{equation*}
for all $x\in X_{\lceil\veps^{-1}\rceil}$.
\item[(H2)] For every $k\in\N_+$ there exist $a_k>0$ and $\kappa(k)\in\N_+$ such that
\begin{equation*}
|\la x,D(v)\ra|\leq a_k\cdot\lno x\rno\cdot q_{\kappa(k)}(v), 
\end{equation*}
for all $x \in X_k$ and $v \in E$.
\item[(H3)] For every $k\in\N_+$ there exists $b_k>0$ such that
\begin{equation*}
\lno x\rno\leq b_k\cdot\sup\{ |\la x,D(v)\ra| \mid v\in E\text{ and }q_k(v)\leq1\}, 
\end{equation*}
for all $x \in X_k$.
\item[(H4)] For every $k\in\N_+$ and $x\in X\setminus X_k$ there is $v\in E$ such that $q_k(v)=0$ and $\la x,D(v)\ra\neq0$.
\end{enumerate}
For $F\in X[\calG_{\calC}]^*$ define
\begin{equation*}
\vvvert F\vvvert_k=\sup\{ |\la x,F\ra| \mid x\in X_k\text{ and }\lno x\rno\leq1\}.
\end{equation*}
Then $X[\calG_{\calC}]^*$, equipped with its strong topology, is a Fr\'echet space whose topology is generated by $\la\vvvert\cdot\vvvert_k\ra_k$, and $D:E\to X[\calG_{\calC}]^*$ is continuous and surjective.
Furthermore, given $F\in X[\calG_{\calC}]^*$ and $k\in\N_+$ such that $\vvvert F\vvvert_k>0$, for every $\veps>0$ there exists $v\in E$ such that
\begin{equation*}
D(v)=F \qquad\text{and}\qquad q_k(v)\leq (1+\veps)\cdot b_k\cdot\vvvert F\vvvert_k.
\end{equation*}
One may even choose such solutions continuously in $F$ on the set where $\vvvert F\vvvert_k>0$, according to \cite[Theorem 8.1(C4)]{DEP.26c}.
Notice that $\kappa(k)$ occurs in (H2), which ensures the continuity of $D$, but not in the quantitative conclusion, which uses (H3).
\par
{\it Sketch of proof.}
The compactness of the $C_k$'s implies the semireflexivity of $X[\calG_{\calC}]$.
The range of $D$ is dense by the Hahn--Banach theorem, using (H3), and closed by the closed range theorem, using (H4) and the Krein--\v{S}mulian theorem to reduce to compactness of the $C_k$'s.
The quantitative conclusion follows from a further Hahn--Banach separation argument, see \cite[Section 8]{DEP.26c}.
\par
We also recall the notion of uniform sequential density and the extension theorem \cite[5.6--5.7]{DEP.26c}.
A linear subspace $Y\subset X$ is {\em uniformly sequentially $\calG_{\calC}$-dense} if there exist a map $\gamma:\N_+\to\N_+$ and a constant $\Lambda>0$ such that every $x\in X_k$ admits a sequence $\la y_j\ra_j$ in $Y\cap X_{\gamma(k)}$ with
\begin{equation*}
\lim_j\bG(x-y_j)=0 \qquad\text{and}\qquad \sup_j\lno y_j\rno\leq\Lambda\cdot\lno x\rno.
\end{equation*}
If $Y$ has this property and a linear form $F:Y\to\R$ satisfies, for each $\veps>0$, an inequality
\begin{equation*}
|F(y)|\leq\theta\cdot\bG(y)+\veps\cdot\lno y\rno, 
\end{equation*}
$y\in Y\cap X_{\lceil\veps^{-1}\rceil}$, for some $\theta>0$, then $F$ admits a unique $\calG_{\calC}$-continuous linear extension $\hat F:X\to\R$.
More precisely, with $\ell=\max\{\lceil\Lambda k\rceil,\gamma(k)\}$, one has
\begin{equation*}
\sup_{x\in C_k}|\hat F(x)|\leq\sup_{y\in Y\cap C_\ell}|F(y)|.
\end{equation*}
In the particular case $\Lambda=1$ and $\gamma(k)=k+1$, each $x\in X_k$ with $\lno x\rno\leq1$ can be approximated by $y_j\in Y\cap X_{k+1}$ with $\lno y_j\rno\leq1$, so that
\begin{equation*}
\vvvert\hat F\vvvert_k\leq\sup\{|F(y)| \mid y\in Y\cap X_{k+1}\text{ and }\lno y\rno\leq1\}.
\end{equation*}
This is the quantitative extension estimate used later in \ref{9.4}.
\end{Empty}

\begin{Empty}[Newtonian potential]
\label{Newton}
We let $\Gamma$ denote the Newtonian kernel in $\Rm$, normalised so that $\triangle\Gamma=\bdelta_0$.
In particular, $|\nabla\Gamma(x)|_2=\frac{1}{m\balpha(m)}\cdot|x|_2^{1-m}$ for $x\neq 0$.
Let $f:\Rm\to\R$ be Lipschitz, $a\in\Rm$, $R>0$, and assume that $\rmspt f\subset B(a,R)$.
Put $g=\Gamma*f$.
Then $g\in C^2(\Rm)$ and $\triangle g=f$.
If, moreover, $f\in\calD(\Rm)$, then $g\in C^\infty(\Rm)$.
\begin{enumerate}
\item[(A)] {\it For every $i,j\in\{1,\ldots,m\}$ the following hold:
\begin{enumerate}
\item[(a)] $\displaystyle \|D^2_{ij}g\|_\infty\leq(3m-2)\cdot R\cdot(\rmLip f)$;
\item[(b)] If $x\in\Rm\setminus B(a,2R)$, then $\displaystyle |D^2_{ij}g(x)|\leq \frac{(m-1)2^m}{m\balpha(m)}\cdot\frac{\|f\|_{L_1}}{|x-a|_2^m}$.
\end{enumerate}
}
\item[(B)] {\it $\displaystyle \|\nabla g\|_\infty\leq R\cdot\|f\|_\infty$.}
\end{enumerate}
\par
{\it Proof.}
First, $\|f\|_\infty\leq R\cdot(\rmLip f)$, since $f$ vanishes on $\rmbdry B(a,R)$.
If $|x-a|_2\geq 2R$ and $y\in\rmspt f$, then $|x-y|_2\geq |x-a|_2/2$.
Since
\begin{equation*}
|D^2_{ij}\Gamma(z)|\leq\frac{m-1}{m\balpha(m)}\cdot|z|_2^{-m}
\end{equation*}
when $z \neq 0$, according to \cite[(2.14) p.17]{GILBARG.TRUDINGER}, differentiation under the integral sign gives (A)(b).
It follows also that, for $|x-a|_2\geq2R$,
\begin{equation*}
|D^2_{ij}g(x)|\leq\frac{m-1}{m}\cdot\|f\|_\infty\leq\frac{m-1}{m}\cdot R\cdot(\rmLip f).
\end{equation*}
For $x\in B(a,2R)$, the representation formula in \cite[Lemma 4.2]{GILBARG.TRUDINGER}, with radius $3R$, gives
\begin{equation*}
|D^2_{ij}g(x)|\leq (m-1)\cdot\int_0^{3R}\frac{\rmosc_f(t)}{t}\,dt+\frac{1}{m}\cdot\|f\|_\infty\leq\left(3m-3+\frac{1}{m}\right)\cdot R\cdot(\rmLip f),
\end{equation*}
where the integral outside $B(x,3R)$ vanishes because $\rmspt f\subset B(x,3R)$.
These two bounds prove (A)(a).
Finally, for each $x\in\Rm$,
\begin{equation*}
|\nabla g(x)|_2\leq\frac{\|f\|_\infty}{m\balpha(m)}\cdot\int_{B(a,R)}|x-y|_2^{1-m}\,d\ssfL^m(y)\leq R\cdot\|f\|_\infty,
\end{equation*}
where the last inequality follows since the integral of a radially decreasing function over a ball is maximal when the ball is centred at the origin.
This proves (B).\cqfd
\end{Empty}

\begin{Empty}[Ornstein non-inequalities]
\label{Ornstein}
We recall the following theorem of Ornstein \cite{ORN.62} and we refer to Kirchheim-Kristensen \cite{KIR.KRI.16} for a different proof and to Faraco-Guerra \cite{FAR.GUE.22} for a simplified argument when $m=2$.
See also Conti-Faraco-Maggi \cite{CON.FAR.MAG.05} for similar constructions.
We let $\rmHom(\rmSym(m),\R^d)$ denote the set of linear homomorphisms $\rmSym(m) \to \R^d$, $d \in \N_+$.
\begin{enumerate}
\item[(A)] {\it 
Let $P_i \in \rmHom(\rmSym(m),\R^{d_i})$, $i=1,2$.
If there exists $\bc > 0$ such that
\begin{equation*}
\int_{\Rm} |P_1(D^2 \vphi)| \,d\ssfL^m \leq \bc \cdot \int_{\Rm} |P_2(D^2(\vphi))| \,d\ssfL^m
\end{equation*}
for all $\vphi \in \calD(\Rm)$ such that $\rmspt \vphi \subset B(0,1)$ then $P_1 = L \circ P_2$ for some $L \in \rmHom(\R^{d_2},\R^{d_1})$.
}
\end{enumerate}
\par 
{\it Proof.}
For the case $m=2$, see \cite[theorem 0.4]{FAR.GUE.22}.
The general case is a corollary of \cite[theorem 1.3]{KIR.KRI.16}.\cqfd
\par 
We define the trace-free Hessian of a test function $\vphi \in \calD(\Rm)$ by the formula $D^2_0\vphi = D^2 \vphi - \frac{1}{m} \cdot (\triangle \vphi) \cdot I_m$ (recall \ref{np.matrices}).
\begin{enumerate}
\item[(B)] {\it $(\forall \bc > 0)(\exists \vphi \in \calD(\Rm)): \rmspt \vphi \subset B(0,1),\, \|D^2_0\vphi\|_{L_1} = 1$ and $\bc < \|\triangle \vphi\|_{L_1}$.
}
\end{enumerate}
\par 
{\it Proof.}
Let $P_1 : \rmSym(m) \to \R$ be $P_1 = \rmtrace$ and $P_2 : \rmSym(m) \to \R^{m \times m}$ be $P_2(A) = A - \frac{1}{m} \cdot (\rmtrace A) \cdot I_m$.
If $\neg (B)$, then $P_1 = L \circ P_2$, by (A), in particular, $\ker P_2 \subset \ker P_1$, contradicting $P_1(I_m) = m$ and $P_2(I_m)=0$.\cqfd
\end{Empty}

\section{Necessary condition for membership in $\rmdiv(C^1(\Rm;\Rm))$}

\begin{Empty}[The norms $\bG$ and $\lno \cdot \rno$ of test functions]
\label{3.1}
Given $u \in \calD_\#(\Rm)$ we define 
\begin{equation*}
\bG(u) = \sup \left\{ \int_{\Rm} u \cdot \vphi \, d\ssfL^m \mid \vphi \in \calD(\Rm) \text{ and } \|\nabla \vphi \|_\infty \leq 1 \right\}.
\end{equation*}
Notice that $\bG$ is a norm on $\calD_\#(\Rm)$.
Next, we put
\begin{equation*}
\lno u \rno = \sum_{i=1}^m \bG(D_i u)
\end{equation*}
which makes sense, since each $D_i u \in \calD_\#(\Rm)$.
This also defines a norm on $\calD_\#(\Rm)$, as each member of $\calD_\#(\Rm)$ whose gradient vanishes is null.
\begin{enumerate}
\item[(A)] {\it For each $u \in \calD_\#(\Rm)$ we have
\begin{equation*}
\lno u \rno = \sup \left\{ \int_{\Rm} u \cdot \rmdiv v \, d\ssfL^m \mid v \in C_c^1(\Rm;\Rm) \text{ and } \max_i \|\nabla v_i \|_\infty \leq 1 \right\}.
\end{equation*}
}
\end{enumerate}
\par 
{\it Proof.}
For every $v \in \calD(\Rm;\Rm)$ we have
\begin{equation*}
\int_{\Rm} u \cdot \rmdiv v \, d\ssfL^m = - \sum_{i=1}^m \int_{\Rm} (D_i u) \cdot v_i \, d\ssfL^m
\end{equation*}
and each term of this sum is optimised over $v_i$ independently of the others.
\par
Now let $v \in C_c^1(\Rm;\Rm)$ and put $v_\delta = \Phi_\delta * v$ for $\delta > 0$.
Then $v_\delta \in \calD(\Rm;\Rm)$ and $\max_i\|\nabla (v_\delta \ip e_i)\|_\infty \leq \max_i\|\nabla v_i\|_\infty$.
Moreover, $\rmdiv v_\delta \to \rmdiv v$ uniformly as $\delta \to 0^+$.
Applying the preceding equality to $v_\delta$ and passing to the limit shows that the supremum is unchanged when the smooth vector fields are replaced by those in $C_c^1(\Rm;\Rm)$.\cqfd
\begin{enumerate}
\item[(B)] {\it For each $u \in \calD_\#(\Rm)$ we have
\begin{equation*}
\bG(u) \leq \rmdiam(\rmspt u) \cdot \|u\|_{L_1} \quad\text{and }\quad \lno u \rno \leq m \cdot \|u\|_{L_1}.
\end{equation*}
}
\end{enumerate}
\par 
{\it Proof.}
{\bf (i)}
Let $\vphi \in \calD(\Rm)$.
Assuming without loss of generality that $u \neq 0$, pick $x_0 \in \rmspt u$.
As $|\vphi(x)-\vphi(x_0)| \leq |x-x_0|_2 \cdot \|\nabla \vphi\|_\infty$, for all $x \in \Rm$,
and $\int_{\Rm} u \, d\ssfL^m = 0$, we have
\begin{equation*}
\int_{\Rm} \vphi \cdot u \, d\ssfL^m = \int_{\Rm} (\vphi(x)-\vphi(x_0))\cdot u(x) \, d\ssfL^m(x) \leq \rmdiam(\rmspt u) \cdot \|\nabla \vphi \|_\infty \cdot \|u\|_{L_1}.
\end{equation*}
This proves the first inequality.
\par 
{\bf (ii)}
Let $v \in \calD(\Rm;\Rm)$.
As
\begin{equation*}
\int_{\Rm} u \cdot \rmdiv v \, d\ssfL^m = \sum_{i=1}^m \int_{\Rm} u \cdot (D_i v_i) \, d\ssfL^m \leq \sum_{i=1}^m \|u\|_{L_1} \cdot \| D_i v_i \|_\infty,
\end{equation*}
the second inequality is a consequence of (A).\cqfd
\end{Empty}

\begin{Lemma}
\label{3.2}
There are constants $\bc_{\theTheorem}(k) > 0$ with the following property.
For every $v \in C^1(\Rm;\Rm)$ and $k \in \N_+$ there exists $\hat{v} \in C_c^1(\Rm;\Rm)$ such that $\rmspt \hat{v} \subset U(0,k+1)$, $\rmdiv \hat{v} = \rmdiv v$ in a neighbourhood of $B(0,k)$, and
\begin{equation*}
\max_i \|\nabla \hat{v}_i\|_\infty
 \leq \bc_{\theTheorem}(k) \cdot \bq_{k+1}(v).
\end{equation*}
\end{Lemma}

\begin{proof}
Let $\psi_k$ be as in \ref{np.omega}.
Define $\hat{v} = \psi_k \cdot (v - v(0))$. 
The first two claims about $\hat{v}$ are readily satisfied.
For each $i \in \{1,\ldots,m\}$ we have $\nabla \hat{v}_i = (v_i - v_i(0)) \cdot \nabla \psi_k + \psi_k \cdot \nabla v_i$.
Since $\rmspt \psi_k \subset B(0,k+1)$ and $|v_i(x)-v_i(0)| \leq |x|_2 \cdot \bq_{k+1}(v)$ for all $x \in B(0,k+1)$, we obtain $\|\nabla \hat{v}_i \|_\infty \leq (k+1) \cdot \bq_{k+1}(v) \cdot \|\nabla \psi_k\|_\infty + \bq_{k+1}(v) \leq (2k+3) \cdot \bq_{k+1}(v)$.
Thus, the lemma holds with $\bc_{\theTheorem}(k) = 2k+3$.
\end{proof}

\begin{Theorem}
\label{3.3}
Let $v \in C^1(\Rm;\Rm)$.
The following holds.
\begin{multline*}
(\forall \veps > 0)(\forall k \in \N_+)(\exists \theta > 0)(\forall u \in \calD_\#(\Rm)): \\ \rmspt u \subset B(0,k) \quad\Rightarrow\quad |\la u , \rmdiv v \ra| \leq \theta \cdot \bG(u) + \veps \cdot \lno u \rno.
\end{multline*}
\end{Theorem}

\begin{proof}
Let $\veps > 0$ and $k \in \N_+$.
Choose $\delta > 0$ such that $\bq_{k+1} (\Phi_\delta * v - v) \leq \veps \cdot \bc_{\ref{3.2}}(k)^{-1}$.
Abbreviate $w = \Phi_\delta * v - v$ and let $\hat{w}$ be associated with $w$ and $k$ in lemma \ref{3.2}.
For every $u \in \calD_\#(\Rm)$ supported in $B(0,k)$ we have
\begin{multline*}
| \la u , \rmdiv v \ra | \leq | \la u , \rmdiv \Phi_\delta * v \ra | + |\la u , \rmdiv (\Phi_\delta*v - v) \ra | = |\la u , \psi_k \cdot \rmdiv \Phi_\delta * v \ra | + |\la u , \rmdiv \hat{w} \ra| \\ 
\leq \| \nabla [ \psi_k \cdot \rmdiv \Phi_\delta * v]\|_\infty \cdot \bG(u) + \left( \max_i \|\nabla \hat{w}_i \|_\infty\right) \cdot \lno u \rno \leq \theta \cdot \bG(u) + \veps \cdot \lno u \rno,
\end{multline*}
where $\theta = \| \nabla [ \psi_k \cdot \rmdiv \Phi_\delta * v]\|_\infty$.
\end{proof}

\section{Preparatory estimates}

\begin{Empty}[$\bG$ and $\lno \cdot \rno$ for distributions]
\label{4.2}
Given a distribution $T \in \calD_0(\Rm)$ with compact support such that $\chi(T) = 0$ we define
\begin{equation*}
\bG(T) = \sup \left\{ \la \vphi , T \ra \mid \vphi \in \calD(\Rm) \text{ and } \|\nabla \vphi \|_\infty \leq 1 \right\} \in [0,\infty]
\end{equation*}
and
\begin{equation*}
\lno T \rno = \sum_{i=1}^m \bG(D_iT) \in [0,\infty].
\end{equation*}
For each $k \in \N_+$ we define vector spaces
\begin{equation*}
X_k = \calD_0(\Rm) \cap \left\{ T \mid \rmspt T \subset B(0,k) , \, \chi(T) = 0 , \text{ and } \lno T \rno < \infty \right\} \quad\text{and}\quad X = \bigcup_{k \in \N_+} X_k.
\end{equation*}
\begin{enumerate}
\item[(A)] {\it For each $T \in X$ we have
\begin{equation*}
\lno T \rno = \sup \left\{ \la \rmdiv v , T \ra  \mid v \in \calD(\Rm;\Rm) \text{ and } \max_i \|\nabla v_i \|_\infty \leq 1 \right\}.
\end{equation*}
Moreover, $\lno \cdot \rno$ is a norm on $X$.
} 
\end{enumerate}
\par 
{\it Proof.}
The argument for smooth vector fields in the proof of \ref{3.1}(A) applies verbatim to prove the first assertion.
Moreover, if $\lno T \rno = 0$ then, upon testing against $v = \vphi \cdot e_i$ for each $\vphi \in \calD(\Rm)$ we infer that $D_iT=0$, for all $i \in \{1,\ldots,m\}$.
Hence, $T$ is constant and, having compact support, it is null.\cqfd
\begin{enumerate}
\item[(B)] {\it If $T \in \calD_0(\Rm)$, $\rmspt(T)$ is compact, $\chi(T)=0$, and $\bM(T) < \infty$ then
\begin{equation*}
\bG(T) \leq \rmdiam(\rmspt T) \cdot \bM(T) \quad\text{and }\quad \lno T \rno \leq m \cdot \bM(T).
\end{equation*} 
}
\end{enumerate}
\par 
{\it Proof.}
Analogous to that of \ref{3.1}(B).\cqfd
\begin{enumerate}
\item[(C)] {\it If $T \in \calD_0(\Rm)$, $\rmspt T$ is compact, $\chi(T)=0$, and $\bG(T) < \infty$ then
\begin{equation*}
\bG(T) = \min \left\{ \bM(S) \mid S \in \calD_1(\Rm) \text{ and } \partial S = T \right\}.
\end{equation*} 
}
\end{enumerate}
\par 
{\it Proof.}
{\bf (i)} If $S \in \calD_1(\Rm)$ and $\partial S = T$ then $\la \vphi , T \ra = \la d\vphi , S \ra \leq \bM(d\vphi) \cdot \bM(S)$ for all $\vphi \in \calD(\Rm)$.
As $\bM(d\vphi) = \|\nabla \vphi\|_\infty$, we have $\bG(T) \leq \bM(S)$.
This shows that $\bG(T)$ is bounded above by the infimum of $\bM(S)$ for such $S$, if any, and it remains to establish the converse inequality and that the infimum is attained.
\par 
{\bf (ii)}
Let $E = \calD^1(\Rm) \cap \{ d\vphi \mid \vphi \in \calD(\Rm)\} \subset C_0(\Rm;\bigwedge^1\Rm)$ and observe that $L : E \to \R$ is linear and well defined by the formula $\la d\vphi , L \ra = \la \vphi , T \ra$: if $d\vphi_1 = d\vphi_2$ then $\vphi_1-\vphi_2$ is constant, hence null, since it is compactly supported.
Furthermore, $|\la d\vphi , L \ra| \leq \bG(T) \cdot \|d\vphi\|_\infty$.
According to Hahn-Banach's theorem, $L$ admits an extension $\hat{L}$ to $C_0(\Rm;\bigwedge^1\Rm)$ such that $|\la \omega , \hat{L} \ra| \leq \bG(T) \cdot \|\omega\|_\infty$.
By the Riesz-Markov representation theorem, there exists an $\bigwedge_1\Rm$-valued Borel measure $\mu$ on $\Rm$, with polar decomposition $\vec{\mu} \cdot \|\mu\|$, such that $\la \omega , \hat{L} \ra = \int_{\Rm} \la \omega, \vec{\mu} \ra \, d\|\mu\|$ for all $\omega \in C_0(\Rm;\bigwedge^1\Rm)$, and $\|\hat{L}\| = \|\mu\|(\Rm) = \bG(T)$.
Let $S = \|\mu\| \wedge \vec{\mu} \in \calD_1(\Rm)$ and note that $\partial S = T$ and $\bM(S) = \|\mu\|(\Rm) = \bG(T)$.\cqfd
\begin{enumerate}
\item[(D)] {\it If $T \in \calD_0(\Rm)$, $\rmspt T$ is compact, $\chi(T)=0$, and $\lno T \rno < \infty$ then
\begin{equation*}
\lno T \rno = \min \left\{ \sum_{i=1}^m \bM(S_i) \mid S_i \in \calD_1(\Rm) \text{ and } \partial S_i = D_iT \text{ for all } i \in \{1,\ldots,m\} \right\}.
\end{equation*} 
}
\end{enumerate}
\par 
{\it Proof.}
For each $i \in \{1,\ldots,m\}$ note that $D_iT$ is a compactly supported distribution with $\chi(D_iT)=0$, and $\bG(D_iT) < \infty$, since $\lno T \rno < \infty$.
Thus, (C) applies to each $D_iT$ and the conclusion follows.\cqfd
\begin{enumerate}
\item[(E)] {\it For $u \in \calD_\#(\Rm)$ let $\iota(u) = \ssfL^m \hel u \in \calD_0(\Rm)$. Then $\iota : \calD_\#(\Rm) \to X$ is linear and injective, and for every $u \in \calD_\#(\Rm)$:
\begin{enumerate}
\item[(a)] $\rmspt \iota(u) = \rmspt u$ and $\chi(\iota(u)) = 0$;
\item[(b)] $\bG(\iota(u)) = \bG(u)$ and $\lno \iota(u) \rno = \lno u \rno$, where the expressions on the right are those defined in \ref{3.1};
\item[(c)] $\iota(u) \in X_k$ if and only if $\rmspt u \subset B(0,k)$, $k \in \N_+$.
\end{enumerate}
}
\end{enumerate}
\par 
{\it Proof.}
Linearity is clear, and if $\iota(u) = 0$ then $\int_{\Rm} u \cdot \vphi \, d\ssfL^m = 0$ for all $\vphi \in \calD(\Rm)$, whence $u = 0$.
(a) is standard and follows from \ref{np.currents}.
Regarding (b), $\la \vphi , \iota(u) \ra = \int_{\Rm} u \cdot \vphi \,d\ssfL^m$ and $D_i \iota(u) = \iota(D_iu)$ for all $\vphi \in \calD(\Rm)$ and $i \in \{1,\ldots,m\}$, so that the definitions of \ref{3.1} and above coincide.
Finally, $\lno \iota(u) \rno = \lno u \rno \leq m \cdot \|u\|_{L_1} < \infty$, by (b) and \ref{3.1}(B), and (c) follows from (a).\cqfd
\par 
Accordingly, we henceforth identify $u$ with $\iota(u)$ whenever convenient, and we regard $\calD_\#(\Rm)$ as a linear subspace of $X$.
In particular, $X_k \cap \calD_\#(\Rm) = \calD_\#(\Rm) \cap \{ u \mid \rmspt u \subset B(0,k)\}$.
\end{Empty}

\begin{Theorem}[Trace-free Hessian estimate]
\label{4.2F}
For every $\vphi \in \calD(\Rm)$ we have $\triangle\vphi \in \calD_\#(\Rm)$ and
\begin{equation*}
\lno \triangle\vphi \rno \leq \frac{m}{m-1} \cdot \|D^2_0 \vphi\|_{L_1}.
\end{equation*}
\end{Theorem}

\begin{proof}
First, notice that $\int_{\Rm} \triangle \vphi \,d\ssfL^m = 0$.
We shall apply \ref{4.2}(D) to $T = \iota(\triangle \vphi) = \ssfL^m \hel \triangle \vphi$.
For each $i \in \{1,\ldots,m\}$ define a smooth compactly supported vector field $w_i = (D^2_0 \vphi) \cdot e_i $ .
Thus,
\begin{equation*}
\rmdiv w_i = \sum_{j=1}^m D_j \left( D_iD_j\vphi - \frac{1}{m} \cdot \delta_{ij}(\triangle \vphi) \right) = \frac{m-1}{m} \cdot D_i(\triangle \vphi).
\end{equation*}
Consequently, if we let $v_i = - \left(\frac{m}{m-1}\right) \cdot w_i$ and $S_i = \ssfL^m \wedge v_i \in \calD_1(\Rm)$ then one easily checks that $\partial S_i = D_i(\ssfL^m \hel \triangle \vphi)$ and $\bM(S_i) = \int_{\Rm} |v_i|_2 \,d\ssfL^m$.
Finally,
\begin{equation*}
\lno \ssfL^m \hel \triangle \vphi \rno \leq \sum_{i=1}^m \bM(S_i) = \left( \frac{m}{m-1} \right) \cdot \| D^2_0 \vphi\|_{L_1},
\end{equation*}
according to \ref{4.2}(D).
\end{proof}

\begin{Empty}[Translations and homotheties]
\label{4.3}
Given $h \in \Rm$ and $r > 0$ we define the translation $\btau_h : \Rm \to \Rm : x \mapsto x + h$ and the homothety $\bmu_r : \Rm \to \Rm : x \mapsto r \cdot x$.
These are readily bijective; $\btau_h^{-1} = \btau_{-h}$ and $\bmu_r^{-1} = \bmu_{r^{-1}}$.
Recall that if $\vphi \in \calD(\Rm)$ then $(\btau_h^\# \vphi)(x) = (\vphi \circ \btau_h)(x)=\vphi(x+h)$ and $(\bmu_r^\#\vphi)(x) = (\vphi \circ \bmu_r)(x) = \vphi(r \cdot x)$, $x \in \Rm$.
If $T \in \calD_0(\Rm)$ then $\btau_{h\#}T$ and $\bmu_{r\#}T$ are defined by $\la \vphi , \btau_{h\#}T \ra = \la \btau_h^\# \vphi , T \ra$ and $\la \vphi , \bmu_{r\#}T \ra = \la \bmu_r^\#\vphi , T \ra$.
In particular, if $T \in \calD_0(\Rm)$ is compactly supported then so are $\btau_{h\#}T$ and $\bmu_{r\#}T$, in fact, $\rmspt (\btau_{h\#}T) = \btau_h(\rmspt T)$ and $\rmspt(\bmu_{r\#}T) = \bmu_r(\rmspt T)$.
If, furthermore, $\chi(T) = 0$ then $\chi(\btau_{h\#}T) = \chi(\bmu_{r\#}T)=0$ as well.
It is obvious that $\bM$, $\bG$, and $\lno \cdot \rno$ are invariant under $\btau_{h\#}$.
Recall also that $\bM(\bmu_{r\#}T) = \bM(T)$.
\begin{enumerate}
\item[(A)] {\it If $T \in X$ and $r > 0$ then $\bG(\bmu_{r\#}T) = r \cdot \bG(T)$.}
\end{enumerate}
{\it Proof.}
Immediate from the fact that $\|\nabla(\vphi \circ \bmu_r)\|_\infty = r \cdot \|\nabla \vphi\|_\infty$ and that $\bmu_r$ is bijective.\cqfd
\begin{enumerate}
\item[(B)] {\it If $T \in X$ and $r > 0$ then $\lno \bmu_{r\#} T \rno = \lno T \rno$.}
\end{enumerate}
{\it Proof.}
Let $i \in \{1,\ldots,m\}$ and $\vphi \in \calD(\Rm)$ and note that $D_i(\vphi \circ \bmu_r) = r \cdot (D_i\vphi) \circ \bmu_r$, hence, $\la \vphi , D_i (\bmu_{r\#}T) \ra = - \la D_i\vphi , \bmu_{r\#} T \ra = - \la (D_i\vphi) \circ \bmu_r , T \ra = - \la r^{-1} \cdot D_i(\vphi \circ \bmu_r) , T \ra = \la \vphi , r^{-1} \cdot \bmu_{r\#}(D_iT) \ra$.
Thus, $\lno \bmu_{r\#} T \rno = \sum_{i=1}^m \bG(D_i(\bmu_{r\#}T)) = \sum_{i=1}^m \bG(r^{-1} \cdot \bmu_{r\#}(D_iT)) = \sum_{i=1}^m \bG(D_iT)=\lno T \rno$, by (A).\cqfd
\end{Empty}

\begin{Lemma}
\label{4.4}
Let $f : \Rm \to \R$ be Lipschitz, $a \in \Rm$, $R > 0$, and assume that $\rmspt f \subset B(a,R)$.
For all $\Lambda > 2R$ there exists $w \in C^1(\Rm;\Rm)$ such that $\rmdiv w = f$ on $U(a,\Lambda)$, $\rmspt w \subset B(a,2\Lambda)$, and 
\begin{equation*}
\|\nabla w_i\|_\infty \leq \sqrt{m} \cdot (3m-2) \cdot R \cdot (\rmLip f) + \frac{2^m \|f\|_{L_1}}{m \balpha(m)} \cdot \Lambda^{-m}
\end{equation*}
for all $i \in \{1,\ldots,m\}$.
Furthermore, if $f \in \calD(\Rm)$ then one may choose $w \in \calD(\Rm;\Rm)$.
\end{Lemma}

\begin{proof}
{\bf (i)}
Define $v = \nabla (\Gamma * f)$.
Then $v \in C^1(\Rm;\Rm)$ and $\rmdiv v = f$ in $\Rm$, by \ref{Newton}, and
\begin{equation}
\label{eq.3}
\|\nabla v_i\|_\infty \leq \sqrt{m} \cdot (3m-2) \cdot R \cdot (\rmLip f),
\end{equation}
according to \ref{Newton}(A)(a).
\par 
{\bf (ii)}
Let $x \in \Rm \setminus B(a,2R)$.
If $y \in \rmspt f$ then $|x-y|_2 \geq |x-a|_2 - |y-a|_2 \geq |x-a|_2 - R \geq \frac{|x-a|_2}{2}$ so that 
\begin{multline*}
|v(x)|_2 = \left| \int_{\Rm} (\nabla \Gamma)(x-y) \cdot f(y) \,d\ssfL^m(y) \right|_2 \leq \int_{\Rm} \frac{1}{m\balpha(m)} \frac{1}{|x-y|_2^{m-1}} \cdot |f(y)| \, d\ssfL^m(y) \\ \leq \frac{1}{m\balpha(m)} \frac{2^{m-1}}{|x-a|_2^{m-1}} \cdot \|f\|_{L_1},
\end{multline*}
by \cite[Lemma 4.1 and (2.13) p.17]{GILBARG.TRUDINGER}.
\par 
{\bf (iii)}
Choose $\psi \in \calD(\Rm)$ such that $0 \leq \psi \leq 1$, $\psi = 1$ in a neighbourhood of $B(a,\Lambda)$, $\rmspt \psi \subset B(a,2\Lambda)$, and $\|\nabla \psi\|_\infty \leq 2 \cdot \Lambda^{-1}$.
Put $w = \psi \cdot v$.
Clearly, $\rmspt w \subset B(a,2\Lambda)$ and it readily follows from {\bf (i)} that $w \in C^1(\Rm;\Rm)$ and $\rmdiv w = f$ on $U(a,\Lambda)$.
Moreover, for each $i \in \{1,\ldots,m\}$ we have $\nabla w_i = \psi \cdot \nabla v_i + v_i \cdot \nabla \psi$.
Thus, $\|\nabla w_i\|_\infty \leq \|\nabla v_i\|_\infty + 2 \cdot \Lambda^{-1} \cdot \sup_{x \in \Rm \setminus B(a,\Lambda)} |v_i(x)|$.
The first term is bounded as in \eqref{eq.3} and the second by means of {\bf (ii)}, since $\Lambda > 2R$.
This gives the claimed upper bound for $\|\nabla w_i\|_\infty$.
\par 
{\bf (iv)}
Finally, we note that if $f \in C^{\infty}(\Rm)$ then $\Gamma * f \in C^\infty(\Rm)$, whence, $v,w \in C^\infty(\Rm;\Rm)$.
\end{proof}

\begin{Theorem}[Poincar\'e-type inequality]
\label{4.5}
There are constants $\bc_{\theTheorem}(m,k) > 0$ with the following property.
If $k \in \N_+$, $T \in \calD_0(\Rm)$, $\rmspt T \subset B(0,k)$, and $\chi(T) = 0$ then 
\begin{equation*}
\bG(T) \leq \bc_{\theTheorem}(m,k) \cdot \lno T \rno.
\end{equation*}
\end{Theorem}

\begin{proof}
Let $\vphi \in \calD(\Rm)$ be such that $\|\nabla \vphi \|_\infty \leq 1$.
Define $\hat{\vphi} = \psi_k \cdot (\vphi - \vphi(0)) \in \calD(\Rm)$ and observe that $\la \vphi , T \ra = \la \hat{\vphi} , T \ra$, since $\rmspt T \subset B(0,k)$ and $\chi(T)=0$.
Moreover, $\|\nabla \hat{\vphi}\|_\infty \leq \|\psi_k\|_\infty \cdot \|\nabla \vphi\|_\infty + \left( \sup_{x \in B(0,k+1)} |\vphi(x)-\vphi(0)| \right) \cdot \|\nabla \psi_k\|_\infty \leq 2k + 3$ and $\rmspt \hat{\vphi} \subset B(0,k+1)$.
Let $\Lambda > 2(k+1)$.
Applying lemma \ref{4.4} to $\hat{\vphi}$, $a=0$, and $R = k+1$ yields $w \in \calD(\Rm;\Rm)$ such that, in particular, $\rmdiv w = \hat{\vphi}$ on $U(a,\Lambda)$, hence on a neighbourhood of $\rmspt T$.
Therefore, $\la \hat{\vphi} , T \ra = \la \rmdiv w , T \ra$.
It then follows from \ref{4.2}(A) and the upper bound from lemma \ref{4.4} that
\begin{multline*}
|\la \vphi , T \ra | = | \la \rmdiv w , T \ra | \leq \lno T \rno \cdot \max_i \|\nabla w_i\|_\infty \\ \leq \lno T \rno \cdot \left( \sqrt{m} \cdot (3m-2) \cdot (k+1) \cdot (2k+3) + \frac{2^m \|\hat{\vphi}\|_{L_1}}{m\balpha(m)} \cdot \Lambda^{-m} \right).
\end{multline*}
Letting $\Lambda \to \infty$, we infer from the arbitrariness of $\vphi$ that the theorem holds with $\bc_{\ref{4.5}}(m,k) = \sqrt{m}\cdot(3m-2)\cdot(k+1) \cdot (2k+3)$.
\end{proof}

\begin{Corollary}
\label{4.6}
If $T \in X$ then $\bG(T) < \infty$.
Furthermore, $\lno \cdot \rno$ restricted to $X$ is lower semicontinuous with respect to, respectively, weak* and $\bG$-convergence.
\end{Corollary}

\begin{proof}
The definition of $X$ and theorem \ref{4.5} immediately yield the first conclusion.
The second follows upon noticing that $\bG$-convergence of a sequence implies its weak* convergence and that $\lno \cdot \rno$ is weakly* lower semicontinuous, according to \ref{4.2}(A).
\end{proof}

\begin{Empty}[Smoothing]
\label{4.7}
In this entire number we fix $k \in \N_+$, $T \in X_k$, and $0 < \delta < 1$.
\begin{enumerate}
\item[(A)] {\it $T * \Phi_\delta \in \calD_\#(\Rm)$ and $\rmspt (T * \Phi_\delta) \subset B(0,k+1)$; in particular, $T * \Phi_\delta \in X_{k+1}$, according to \ref{4.2}(E).}
\end{enumerate}
{\it Proof.}
The smoothness of $T * \Phi_\delta$ and the conclusion about its support are standard.
If $\psi \in \calD(\Rm)$ equals 1 on a neighbourhood of $B(0,k+1)$ then $\psi * \Phi_\delta$ equals 1 on a neighbourhood of $B(0,k)$, hence, $\chi(T * \Phi_\delta) = \la \psi , T * \Phi_\delta \ra = \la \psi * \Phi_\delta , T \ra = \chi(T) = 0$.\cqfd
\begin{enumerate}
\item[(B)] {\it $\|T * \Phi_\delta\|_\infty \leq \bc_{\theTheorem}(m) \cdot \delta^{-m} \cdot \|\nabla \Phi\|_\infty \cdot \lno T \rno$ and $\|\nabla(T * \Phi_\delta)\|_\infty \leq \bc_{\theTheorem}(m) \cdot \delta^{-(m+1)} \cdot \max_i \|\nabla D_i\Phi\|_\infty \cdot \lno T \rno$ for some constant $\bc_{\theTheorem}(m) > 0$.}
\end{enumerate}
{\it Proof.}
Let $\psi \in \calD(\Rm)$ be such that $\rmspt \psi \subset B(0,\delta)$.
Recall that $T * \psi$ is smooth and that for fixed $x \in \Rm$ we have $(T * \psi)(x) = \la \vphi_x , T \ra$, where $\vphi_x(y) = \psi(x-y)$.
Notice that $\vphi_x \in \calD(\Rm)$ and $\rmspt \vphi_x \subset B(x,\delta)$ so that an application of lemma \ref{4.4} to $\vphi_x$, $a=x$, and $R=\delta$ yields, for all $\Lambda > \max \{ k + |x|_2 , 2\delta \}$, some $w \in \calD(\Rm;\Rm)$ such that $\rmdiv w = \vphi_x$ on $U(x,\Lambda)$, in particular, on a neighbourhood of $\rmspt T$.
Thus,
\begin{multline*}
|(T * \psi)(x)| = |\la \vphi_x , T \ra| = |\la \rmdiv w , T \ra| \leq \lno T \rno \cdot \max_i \|\nabla w_i\|_\infty \\ \leq \lno T \rno \cdot \left( \sqrt{m} \cdot (3m-2) \cdot \delta \cdot \|\nabla \psi\|_\infty + \frac{2^m\|\vphi_x\|_{L_1}}{m \balpha(m)}\cdot \Lambda^{-m}\right),
\end{multline*}
by \ref{4.2}(A) and \ref{4.4}.
Letting $\Lambda \to \infty$ and applying the inequality to, respectively, $\psi = \Phi_\delta$ and $\psi = D_i \Phi_\delta$ for $i \in \{1,\ldots,m\}$ we infer that
\begin{multline*}
|(T * \Phi_\delta)(x)| \leq \sqrt{m} \cdot (3m-2) \cdot \delta^{-m} \cdot \|\nabla \Phi\|_\infty \cdot \lno T \rno \\ \text{and }\quad |D_i(T * \Phi_\delta)(x)| \leq \sqrt{m} \cdot (3m-2) \cdot \delta^{-(m+1)} \cdot \|\nabla D_i\Phi\|_\infty \cdot \lno T \rno,
\end{multline*}
since $\|\nabla \Phi_\delta\|_\infty = \delta^{-(m+1)} \cdot \|\nabla \Phi\|_\infty$ and $\|\nabla D_i\Phi_\delta\|_\infty = \delta^{-(m+2)} \cdot \|\nabla D_i\Phi\|_\infty$.
As $x$ is arbitrary, the conclusion holds with $\bc_{\theTheorem}(m) = m \cdot (3m-2)$.\cqfd
\begin{enumerate}
\item[(C)] {\it $\lno T * \Phi_\delta \rno \leq \lno T \rno$.}
\end{enumerate}
{\it Proof.}
Let $v \in \calD(\Rm;\Rm)$ and observe that $\|\nabla(\Phi_\delta * v_i)\|_\infty = \|\Phi_\delta * (\nabla v_i)\|_\infty \leq \|\nabla v_i\|_\infty$, by Jensen's inequality, for each $i \in \{1,\ldots,m\}$.
Thus, $|\la \rmdiv v , T * \Phi_\delta \ra| = | \la \rmdiv ( \Phi_\delta * v) , T \ra| \leq \lno T \rno \cdot \max_i \|\nabla (\Phi_\delta * v_i)\|_\infty \leq \lno T \rno \cdot\max_i \|\nabla v_i \|_\infty$, by \ref{4.2}(A), and the conclusion follows from the arbitrariness of $v$.\cqfd
\begin{enumerate}
\item[(D)] {\it $\bG(T - T*\Phi_\delta) \leq \delta \cdot \lno T \rno$.}
\end{enumerate}
{\it Proof.}
Given $\vphi \in \calD(\Rm)$ we start by observing that, for all $x \in \Rm$ we have
\begin{multline*}
\vphi(x) - (\vphi * \Phi_\delta)(x) = \int_{\Rm} \big(\vphi(x) - \vphi(x+h) \big) \cdot \Phi_\delta(h) \, d\ssfL^m(h) \\ = - \int_{\Rm} \left( \int_0^1 h \ip (\nabla \vphi)(x+t \cdot h) \, dt \right) \Phi_\delta(h) \, d\ssfL^m(h)
\\ = - \sum_{i=1}^m \int_{B(0,\delta) \times [0,1]} h_i \cdot D_i(\vphi \circ \btau_{t\cdot h} ) \, d(\ssfL^1 \otimes \ssfL^m \hel \Phi_\delta)(t,h).
\end{multline*}
It follows from \cite[3.26 and 3.27]{RUDIN} that
\begin{multline*}
\la \vphi , T - T * \Phi_\delta \ra = \la \vphi - \vphi * \Phi_\delta , T \ra \\
= \sum_{i=1}^m \int_{B(0,\delta) \times [0,1]} h_i \cdot \la \vphi \circ \btau_{t \cdot h} , D_i T \ra \, d(\ssfL^1 \otimes \ssfL^m \hel \Phi_\delta)(t,h)
\leq \sum_{i=1}^m \delta \cdot \bG(D_iT)\cdot \|\nabla \vphi\|_\infty,
\end{multline*}
since $\|\nabla( \vphi \circ \btau_{t \cdot h})\|_\infty = \|\nabla \vphi\|_\infty$ for all $t$ and $h$.\cqfd
\end{Empty}

\section{Members of $X$}

In this section, we explore the nature of the distributions that are members of $X$.
We recall that a distribution $T \in \calD_0(\Rm)$ belongs to $X$ if and only if it is compactly supported, $\chi(T) = 0$, and $\lno T \rno < \infty$.

\begin{Empty}[Test functions]
\label{5.1}
If $u \in \calD_\#(\Rm)$ then $\iota(u) = \ssfL^m \hel u \in X$, and $\iota(u) \in X_k$ if and only if $\rmspt u \subset B(0,k)$, according to \ref{4.2}(E).
\end{Empty}

\begin{Empty}[Measures]
\label{5.2}
Let $\mu$ be a compactly supported signed Borel measure on $\Rm$ such that $\mu(\Rm) = 0$ and define $T_\mu \in \calD_0(\Rm)$ by $\la \vphi , T_\mu \ra = \int_{\Rm} \vphi \,d\mu$.
Notice that $\rmspt T_\mu = \rmspt \|\mu\|$, $\chi(T_\mu) = \mu(\Rm) = 0$ and $\lno T_\mu \rno \leq m \cdot \bM(T_\mu) = m \cdot \|\mu\|(\Rm) < \infty$, by \ref{4.2}(B).
Thus, $T_\mu \in X$.
\end{Empty}

\begin{Empty}[Members of $X$ are flat chains]
\label{5.3}
Let $T \in X$.
It follows from theorem \ref{4.5} that $\bG(T) < \infty$ and, in turn, from \ref{4.2}(C) that
\begin{equation*}
\bF(T) = \inf \big\{ \bM(R) + \bM(S) \mid T = R + \partial S, R \in \calD_0(\Rm) \text{ and } S \in \calD_1(\Rm) \big\} \leq \bG(T) < \infty,
\end{equation*}
where the first equality is the definition of the flat norm.
This alone, however, does not imply that $T$ is a flat chain, see \eg second paragraph of \ref{5.6}.
Notwithstanding, we infer from \ref{4.7}(D) that $\bF(T - T * \Phi_\delta) \leq \bG(T - T * \Phi_\delta) \to 0$ as $\delta \to 0$.
Since each $T * \Phi_\delta$ is a 0-dimensional normal current, it follows that $T \in \bF_0(\Rm)$ is flat.
\end{Empty}

\begin{Empty}[$X_k$ and the Arens--Eells space of $B(0,k)$]
\label{5.4}
According to \ref{5.3} and the identification recalled in \ref{np.arens-eells},
\begin{equation*}
X_k\subset\calG(B(0,k))\cong\calF(B(0,k)).
\end{equation*}
\end{Empty}

\begin{Empty}[$X_k \subsetneqq \calF(B(0,k))$]
\label{5.5}
Here, we show that the inclusion in \ref{5.4} is strict.
Let $\la r_j \ra_{j \in \N_+}$ be a strictly decreasing sequence of positive real numbers such that $r_1 \leq 1$ and $\sum_{j \in \N_+} r_j < \infty$.
Let also $e_1 \in \Rm$ be the first vector of the canonical basis of $\Rm$.
For each $j \in \N_+$ we define $S_j = \lseg 0 , r_j \cdot e_1 \rseg \in \bF_1(\Rm)$, $T_j = \partial S_j \in \calG(Q)$ and we note that $\bG(T_j) = \bM(S_j) = r_j$.
Since $(\calG(Q),\bG)$ is a Banach space, we infer that $T = \sum_{j \in \N_+} T_j \in \calG(Q)$.
We shall now show that $\lno T \rno = \infty$.
\par 
Fix $n \in \N_+$.
It is not hard to see that there exists $\psi_0 \in C^\infty(\R)$ with the following properties: $\|\psi'_0\|_\infty \leq 1$, $\psi'_0(t) = 1$ for all $0 \leq t \leq r_{n+1}$, $\psi'_0(r_j)=-1$ for all $j \in \{1,\ldots,n\}$, and $\|\psi_0\|_\infty \leq 1$.
Choosing $\psi_1 \in \calD(\Rm)$ such that $0 \leq \psi_1 \leq 1$, $\psi_1$ equals 1 in a neighbourhood of $B(0,1)$ and $\|\nabla \psi_1\|_\infty \leq 1$, we define $\vphi \in \calD(\Rm)$ by $\vphi(x) = \psi_1(x) \cdot \psi_0(x \ip e_1)$.
Note that $\|\nabla \vphi \|_\infty \leq 2$ and $(\nabla \vphi)(x) = \psi_0'(x \ip e_1) \cdot e_1$ for all $x \in B(0,1)$.
Moreover, for each $j \in \N_+$ we have $\la D_1 \vphi , \partial S_j \ra = (D_1\vphi)(r_j \cdot e_1) - (D_1\vphi)(0) = \psi'_0(r_j) - 1$.
This vanishes if $j > n$ and equals $-2$ if $j \leq n$.
Accordingly,
\begin{equation*}
2 \cdot \bG(D_1 T) \geq \la \vphi , D_1T\ra = - \la D_1\vphi , T \ra = - \sum_{j=1}^n \big( \psi'_0(r_j) - \psi'_0(0) \big) = 2n.
\end{equation*}
Thus, $\bG(D_1T) \geq n$. 
Since $n$ is arbitrary, $\lno T \rno \geq \bG(D_1T) = \infty$.
\end{Empty}

\begin{Empty}[Regularity of $S$ and $S_i$'s in \ref{4.2}(C) and (D)]
\label{5.6}
Let $T \in X_k$.
It follows from \ref{5.4} that $T \in \calG(B(0,k))$.
We then infer from \cite[theorem 5.2(C)]{DEP.26b} that, in \ref{4.2}(C) above $S$ can be chosen to be a 1-dimensional flat current.
In other words, there exists $S \in \bF_1(\Rm)$ such that $\partial S = T$, $\bM(S) = \bG(T)$, and $\rmspt S \subset B(0,k)$.
The question then arises whether, in \ref{4.2}(D) above one can choose the $S_i$'s to be 1-dimensional flat currents as well: do there exist $S_i \in \bF_1(\Rm)$ such that $\partial S_i = D_i T$, $i \in \{1,\ldots,m\}$, and $\lno T \rno = \sum_{i=1}^m \bM(S_i)$?
This fails for a simple reason: if $S_i$ were flat then so would $D_iT$ be, but this is not the case for simple choices of $T \in X$ as we show below.
\par 
For instance, let $a,b \in \Rm$ be distinct and consider $T = \bdelta_b - \bdelta_a$.
Then $T \in X$, according to $\ref{5.2}$ and, indeed, $T \in \bN_0(\Rm) \subset \bF_0(\Rm)$.
It also follows from theorem \ref{4.5} or directly from the definition of $\bG$ that $\bF(T) \leq \bG(T) < \infty$.
We now show in two different ways that $D_1T \not \in \bF_0(\Rm)$.
Both arguments proceed by contradiction.
If $D_1T$ were flat then so would $(D_1T) \hel \psi$ be for any $\psi \in \calD(\Rm)$ and choosing $\psi$ equal to 1 in a neighbourhood of $b$ and so that $a \not \in \rmspt \psi$, it would follow that $D_1 \bdelta_b$ is flat as well.
Assume without loss of generality that $b=0$;
Consider $f : \Rm \to \Rm : (x_1,\ldots,x_m) \mapsto (x_1^+,x_2,\ldots,x_m)$ and notice that $f$ is Lipschitz and $f = \rmid_{\Rm}$ on $\rmspt D_1\bdelta_0 = \{0\}$.
Abbreviate $R = D_1\bdelta_0$. If $R$ were flat then $R = \rmid_{\Rm\#}R = f_\#R$, see \cite[4.1.15 end on p. 372]{GMT}.
Moreover, according to \cite[4.1.14]{GMT}, $f_\#R$ would be the (flat, hence, weakly*) limit of $f_{\veps\#}R$ as $\veps \to 0$ where $f_\veps(x_1,\ldots,x_m) = (g_\veps(x_1),x_2,\ldots,x_m)$ and $g_\veps = (\cdot - \veps)^+ *\Phi_{\veps/2}$.
As $\la D_1 f_{\veps}^\# \vphi , \bdelta_0 \ra = D_1(\vphi \circ f_\veps)(0) = 0$ for all $\vphi \in \calD(\Rm)$, we would have $D_1\bdelta_0 = 0$, which is clearly not the case.
\par 
The second proof we offer that $D_1T$ is not flat proceeds via the identification of $\calG(Q)$ with $\calF(Q)$, \cite[7.1]{DEP.26b} or \cite[6.1.1(D)]{DEP.26a}.
If $D_1T$ were flat, we would have $\la u_n , D_1\bdelta_b - D_1\bdelta_a \ra \to 0$ as $n \to \infty$ for every sequence $\la u_n \ra_{n \in \N_+}$ converging pointwise to 0 and such that $\sup_n \rmLip u_n < \infty$.
Assume without loss of generality that $b=0$ and $a_1 \neq 0$, and consider $u_n(x) = n^{-1} \cdot \psi(x) \cdot \sin(n\cdot(x_1-a_1))$, where $\psi \in \calD(\Rm)$ vanishes in a neighbourhood of $b$ and equals $1$ in a neighbourhood of $a$. We would then have $\la D_1u_n , \bdelta_b - \bdelta_a \ra \to 0$ as $n \to \infty$, yet $\la D_1u_n , \bdelta_b - \bdelta_a \ra =-1$ for all $n$.
\end{Empty}

\begin{Lemma}
\label{5.7}
For every closed ball $B \subset \Rm$ and $M > 0$ there exists $u \in \calD_\#(\Rm)$ such that $\rmspt u \subset B$, $\lno u \rno = 1$, and $\|u\|_{L_1} > M$.
\end{Lemma}

\begin{proof}
By \ref{Ornstein}(B), there exists $\vphi \in \calD(\Rm)$ such that $\rmspt \vphi \subset B(0,1)$, $\| D^2_0 \vphi \|_{L_1} = \frac{m-1}{m}$ and $\|\triangle \vphi\|_{L_1} > M$.
Letting $u_0 = \triangle \vphi$, we infer from Theorem \ref{4.2F} that $u_0 \in \calD_\#(\Rm)$, $\lno u_0 \rno \leq 1$, and $\|u_0\|_{L_1} > M$.
Multiplying $u_0$ by a scalar if necessary, we may achieve $\lno u_0 \rno = 1$.
Now, writing $B = B(a,r)$ and putting $T_0 = \ssfL^m \hel u_0$, we see that $T = \btau_{a\#}\bmu_{r\#}T_0$ satisfies $\lno T \rno = \lno T_0 \rno = 1$, by \ref{4.3}(B), and $\rmspt T \subset B$.
Moreover, $T = \ssfL^m \hel u$ where $u(x) = \frac{1}{r^m} \cdot u_0 \left( \frac{x-a}{r}\right)$, $x \in \Rm$, and $\|u\|_{L_1} = \|u_0\|_{L_1} > M$.
\end{proof}

\begin{Empty}[Non-measures in $X$]
\label{5.8}
Here we show that not all distributions in $X$ are associated with a measure as in \ref{5.2}.
Let $\la B_j \ra_j$ be a sequence of pairwise disjoint closed balls contained in $B(0,1)$.
For each $j \in \N_+$ there exists a smooth distribution $T_j = \ssfL^m \hel u_j$ supported in $B_j$ so that $\lno T_j \rno = 1$ and $\bM(T_j) = \|u_j\|_{L_1} > 2^j$, according to lemma \ref{5.7}.
We observe that the sequence of distributions $\left\la \sum_{j=1}^n 2^{-j} \cdot T_j \right\ra_{n \in \N}$ is weakly* convergent, since 
\begin{equation*}
\left| \left\la \vphi , \sum_{j=1}^n 2^{-j} \cdot T_j \right\ra \right| \leq \|\nabla \vphi\|_\infty \sum_{j=1}^n 2^{-j} \cdot \bG(T_j) \leq \|\nabla \vphi \|_\infty \cdot \bc_{\ref{4.5}}(m,1) \sum_{j=1}^n 2^{-j},
\end{equation*}
according to the definition of $\bG$ and theorem \ref{4.5}.
Thus, $T = \sum_{j=1}^\infty 2^{-j} \cdot T_j$ is a distribution, by the Banach--Steinhaus theorem, since $\calD(\Rm)$ is barrelled, see \cite[6.17]{RUDIN}.
Moreover, $\rmspt T \subset B(0,1)$, $\chi(T) = 0$, and $\lno T \rno \leq 1$, by corollary \ref{4.6}.
Thus, $T \in X$.
Assume if possible that there exists a signed Borel measure $\mu$ such that $\la \vphi , T \ra = \int_{\Rm} \vphi \,d\mu$ for all $\vphi \in \calD(\Rm)$.
It is then obvious that $\mu \hel (\rmint B_j) = 2^{-j} \cdot \ssfL^m \hel u_j$ for all $j$, whence, $|\mu|(\Rm) \geq \sum_{j=1} 2^{-j} \cdot \|u_j\|_{L_1} = \infty$, a contradiction.
\end{Empty}

\section{Sufficient condition for membership in $\rmdiv(C^1(\Rm;\Rm))$}

\begin{Empty}[The topology $\calG_{\lno \cdot \rno}$]
\label{6.1}
We recall that $\la X_k \ra_{k \in \N_+}$, $X$, $\bG$, and $\lno \cdot \rno$ have been defined in \ref{4.2}.
The absolutely convex sets
\begin{equation*}
C_k = X_k \cap \{ T \mid \lno T \rno \leq k \}
\end{equation*} 
readily form a localising family $\calC = \{ C_k \mid k \in \N_+\}$, see \cite{DEP.26c}.
\par 
We observe that $\bG$ is a norm on $X$.
First, note that $\bG(T) < \infty$ whenever $T \in X$, by theorem \ref{4.5}.
All the remaining properties of a norm are immediate consequences of the definition of $\bG$.
We let $\calG$ denote the locally convex topology on $X$ associated with $\bG$.
We write $\calG_{\lno \cdot \rno} = \calG_{\calC}$ for the associated localised topology.
\end{Empty}

\begin{Theorem}[Compactness]
\label{6.2}
For every $k \in \N_+$ the space $C_k[\calG \hel C_k]$ is compact, \ie $C_k$ is $\bG$-compact.
\end{Theorem}

\begin{proof}
{\bf (i)} 
We start by showing that $C_k$ is totally bounded with respect to $\bG$.
Let $0 < \veps < 1$ and put $\delta = \frac{\veps}{3k}$.
Consider the set of smooth functions $\calH = C(\Rm) \cap \{ T*\Phi_\delta \mid T \in C_k \}$.
Each member of $\calH$ is supported in $B(0,k+1)$, by \ref{4.7}(A), and $\calH$ is $\|\cdot\|_\infty$-bounded and equi-Lipschitz, by \ref{4.7}(B).
It follows from Ascoli's theorem that there are $T_1,\ldots,T_N \in C_k$ such that for each $T \in C_k$ there exists $n \in \{1,\ldots,N\}$ with $\|T * \Phi_\delta - T_n * \Phi_\delta\|_\infty \leq \frac{\veps}{6 \balpha(m)(k+1)^{m+1}}$.
Therefore, $\bM(T * \Phi_\delta - T_n *\Phi_\delta) = \| T * \Phi_\delta - T_n *\Phi_\delta \|_{L_1} \leq \balpha(m)(k+1)^m \cdot \|T * \Phi_\delta - T_n *\Phi_\delta\|_\infty \leq \frac{\veps}{6(k+1)}$, whence, $\bG(T * \Phi_\delta - T_n *\Phi_\delta) \leq 2(k+1) \cdot \bM(T * \Phi_\delta - T_n *\Phi_\delta) \leq \frac{\veps}{3}$, by \ref{3.1}(B) and \ref{4.2}(E).
Finally, $\bG(T-T_n) \leq \bG(T - T*\Phi_\delta) + \bG(T*\Phi_\delta - T_n * \Phi_\delta) + \bG(T_n *\Phi_\delta - T_n) \leq \delta \cdot \lno T \rno + \frac{\veps}{3} + \delta \cdot\lno T_n \rno \leq \veps$, by \ref{4.7}(D).
\par 
{\bf (ii)}
We now show that $C_k$ is $\bG$-complete.
Let $\la T_j \ra_j$ be $\bG$-Cauchy in $C_k$.
It is readily weakly* Cauchy in $\calD_0(\Rm)$, hence, converges weakly* to some $T \in \calD_0(\Rm)$, by the Banach--Steinhaus theorem, since $\calD(\Rm)$ is barrelled, see \eg \cite[6.17]{RUDIN}.
One readily checks that $\rmspt T \subset B(0,k)$ and $\chi(T) = 0$.
It then follows from \ref{4.6} that $\lno T \rno \leq k$, thus, $T \in C_k$.
Finally, given $\veps > 0$ let $j_0 \in \N_+$ be such that $\bG(T_{j_1}-T_{j_2}) \leq \veps$ whenever $\min\{j_1,j_2\} \geq j_0$.
Fix $j_2 \geq j_0$.
If $\vphi \in \calD(\Rm)$ and $\|\nabla \vphi\|_\infty \leq 1$ then $|\la \vphi , T-T_{j_1}\ra| \leq \veps$, for some $j_1 \geq j_0$, whence, $|\la \vphi , T-T_{j_2}\ra| \leq |\la \vphi , T-T_{j_1}\ra| + \bG(T_{j_1}-T_{j_2}) \leq 2 \veps$.
As $\vphi$ is arbitrary, $\bG(T-T_{j_2}) \leq 2 \veps$ and the proof is complete.
\end{proof}

\begin{Remark}
\label{6.3}
The above theorem also results from the characterisation of $\bG$-compactness in \cite[theorem 8.3]{DEP.26b}.
Specifically, if $T \in \calG(Q)$ and $\veps > 0$, define (see \cite[8.2(D)]{DEP.26b})
\begin{equation*}
\bkappa(T,\veps) = \inf \big\{ \bM(\hat{T}) \mid \hat{T} \in \calG(Q) \text{ and } \bG(T - \hat{T}) < \veps  \big\}.
\end{equation*}
According to \cite[theorem 8.3]{DEP.26b}, a $\bG$-closed set $\calS \subset \calG(Q)$ is $\bG$-compact if and only if it is $\bG$-bounded and there exists a non-increasing $\kappa : \R_{>0} \to \R_{\geq 0}$ such that $\bkappa(T,\veps) \leq \kappa(\veps)$ for all $T \in \calS$ and all $\veps > 0$.
We shall apply this to $\calS_k = \bmu_{r_k\#}(C_k)$, where $r_k = \frac{1}{2k}$.
It is obvious that $\calS_k \subset \calG(Q)$.
Moreover, $\calS_k$ is $\bG$-closed and $\bG$-bounded, by \ref{4.3}(B), since so is $C_k$.
Let $\veps > 0$ and set $\delta = \min \left\{ \frac{1}{3} , \frac{\veps}{2k}\right\}$.
If $T \in \calS_k$, put $\hat{T} = T * \Phi_\delta$ and note that $\hat{T} \in \calG(Q)$ and $\bG(T - \hat{T}) \leq \delta \cdot \lno T \rno < \veps$, by \ref{4.7}(D) and \ref{4.3}(B).
Furthermore, $\bM(\hat{T}) \leq \balpha(m) \cdot \|T * \Phi_\delta\|_\infty \leq \balpha(m) \cdot \bc_{\ref{4.7}}(m) \cdot \|\nabla \Phi\|_\infty \cdot k \cdot \min \left\{ \frac{1}{3} , \frac{\veps}{2k}\right\}^{-m}$, in view of \ref{4.7}(B).
The right member of the inequality defines $\kappa(\veps)$, independently of $T$.
\end{Remark}

\begin{Empty}[The space $\bE^1(\Rm)$]
\label{6.4}
Abbreviate 
\begin{equation*}
\bE^1(\Rm) = C^1(\Rm;\Rm) \cap \{ v \mid v(0) = 0 \}.
\end{equation*}
We consider the sequence of seminorms $\la \bq_k \ra_k$ on $\bE^1(\Rm)$, defined in \ref{np.omega}.
It is separating and defines a metrisable locally convex topology $\calT$.
We observe that this yields a structure of a Fr\'echet space on $\bE^1(\Rm)$.
Indeed, for every $v \in \bE^1(\Rm)$ and $x \in B(0,k)$, $|v(x)|_2 \leq \sqrt{m} \cdot k \cdot \bq_k(v)$.
Thus, if $\la v_j \ra_j$ is $\bq_k$-Cauchy then $\la v_j|_{U(0,k)} \ra_j$ is $\|\cdot\|_\infty$-Cauchy in $C(U(0,k);\Rm)$ and converges uniformly and with respect to $\bq_k$ to a $C^1$ vector field $w_k : U(0,k) \to \Rm$.
As $w_{k+1}|_{U(0,k)}=w_k$, this defines the $\calT$-limit of $\la v_j \ra_j$.
\end{Empty}

\begin{Empty}[Uniform sequential density]
\label{6.5}
Recall from \cite[5.6]{DEP.26c} that a linear subspace $Y \subset X$ is {\em uniformly sequentially $\calG_{\lno \cdot \rno}$-dense} if (by definition) there are $\gamma \in \N^{\N}$ and $\Gamma > 0$ such that for every $k$ and every $T \in X_k$ there exists a sequence $\la T_j \ra_j$ in $Y \cap X_{\gamma(k)}$ with $\lim_j \bG(T-T_j) = 0$ and $\sup_j \lno T_j \rno \leq \Gamma \cdot \lno T \rno$.
Here, we observe that $\calD_\#(\Rm)$ (regarded as a subspace of $X$ by means of $\iota$, see \ref{4.2}(E)) is uniformly sequentially $\calG_{\lno \cdot \rno}$-dense in $X$, with $\gamma(k)=k+1$ and $\Gamma = 1$.
Letting $T_j = T * \Phi_{1/j}$, we note that $T_j \in \calD_\#(\Rm) \cap X_{k+1}$, by \ref{4.7}(A).
Moreover, $\lno T_j \rno \leq \lno T \rno$, by \ref{4.7}(C), and $\bG(T-T_j) \leq \frac{1}{j} \cdot \lno T \rno \to 0$ as $j \to \infty$, according to \ref{4.7}(D).
\end{Empty}

\begin{Empty}[The divergence operator $\bD$]
\label{6.6}
Letting $\calD_\#(\Rm)^\star$ denote the algebraic dual of $\calD_\#(\Rm)$, we consider the linear map
\begin{equation}
\label{eq.4}
\rmdiv : \bE^1(\Rm) \to \calD_\#(\Rm)^\star : v \mapsto \rmdiv v
\end{equation} 
defined by $\la u , \rmdiv v \ra = \int_{\Rm} u \cdot \rmdiv v \, d\ssfL^m$.
\par 
Let $k \in \N_+$ and $T \in X_k$.
We say that a sequence $\la u_j \ra_j$ in $\calD_\#(\Rm)$ is {\em admissible for $T$ at level $k$} if $\rmspt u_j \subset B(0,k+1)$ for all $j$, $\lim_j \bG(T - \iota(u_j)) = 0$, and $\sup_j \lno u_j \rno < \infty$.
Given $\Gamma > 0$, we say that it is {\em $\Gamma$-admissible for $T$ at level $k$} if, moreover, $\sup_j \lno u_j \rno \leq \Gamma \cdot \lno T \rno$.
According to \ref{6.5}, every $T \in X_k$ admits a sequence that is $1$-admissible for $T$ at level $k$.
\par 
{\it 
For every $v \in \bE^1(\Rm)$, every $k \in \N_+$, every $T \in X_k$, and every sequence $\la u_j \ra_j$ admissible for $T$ at level $k$, the limit $\lim_j \la u_j , \rmdiv v \ra$ exists and is independent of the choice of such $\la u_j \ra_j$, and of $k$.
}
\par 
{\it Proof.}
Given $\veps > 0$ and $\theta > 0$ associated with $\veps$ and $k+1$ in theorem \ref{3.3} we infer that $|\la u_{j_1} - u_{j_2} , \rmdiv v \ra| \leq \theta \cdot \bG(\iota(u_{j_1} - u_{j_2})) + 2\cdot \veps \cdot \Gamma$ for all $j_1,j_2 \in \N_+$, where $\Gamma = \sup_j \lno u_j \rno$.
Accordingly, the sequence $\la \la u_j , \rmdiv v \ra \ra_j$ is Cauchy in $\R$.
If $\la \tilde{u}_j \ra_j$ is another sequence admissible for $T$ at level $k$ then by interlacing $\la u_j \ra_j$ and $\la \tilde{u}_j \ra_j$ we obtain a sequence admissible for $T$ at level $k$, whence, $\lim_j \la \tilde{u}_j , \rmdiv v \ra = \lim_j \la u_j , \rmdiv v \ra$.
Finally, if $T \in X_k \subset X_{k'}$, $k \leq k'$, then every sequence admissible for $T$ at level $k$ is admissible for $T$ at level $k'$, and the independence of $k$ follows from the preceding sentence applied at level $k'$.\cqfd
\par 
We are now in a position to define 
\begin{equation*}
\bD : \bE^1(\Rm) \to X^\star
\end{equation*}
in the following way. 
Given $v \in \bE^1(\Rm)$ and $T \in X$ we choose $k$ such that $T \in X_k$ and a sequence $\la u_j \ra_j$ admissible for $T$ at level $k$, and we put
\begin{equation*}
\la T , \bD(v) \ra = \lim_j \la u_j , \rmdiv v \ra.
\end{equation*}
Let $T, \tilde{T} \in X$ and $t \in \R$.
If $\la u_j \ra_j$ and $\la \tilde{u}_j \ra_j$ are admissible for, respectively, $T$ and $\tilde{T}$ at a common level $k$, then $\la t \cdot u_j + \tilde{u}_j \ra_j$ is admissible for $t \cdot T + \tilde{T}$ at level $k$ and it follows from the linearity of $\rmdiv v$ that $\la t \cdot T + \tilde{T} , \bD(v) \ra = t \cdot \la T , \bD(v) \ra + \la \tilde{T} , \bD(v)\ra$, \ie $\bD(v)$ is a member of $X^\star$.
\par 
We note that if $u \in \calD_\#(\Rm)$ and $\rmspt u \subset B(0,k)$ then the constant sequence $u_j = u$ is $1$-admissible for $\iota(u)$ at level $k$, whence, $\la \iota(u) , \bD(v) \ra = \int_{\Rm} u \cdot \rmdiv v \, d\ssfL^m$ for all $v \in \bE^1(\Rm)$.
\end{Empty}

\begin{Empty}
\label{6.7}

Recall that $\bE^1(\Rm)$ is defined in \ref{6.4}, the $\bq_k$'s are defined in \ref{np.omega}, $X$ and $X_k$'s as well as $\bG$ and $\lno \cdot \rno$ are defined in \ref{4.2}, and $\bD$ is defined in \ref{6.6}.
We are now ready to apply the results of \cite{DEP.26c} to the solvability of the equation $\rmdiv v = T$.
Specifically, we shall apply \cite[theorem 8.1]{DEP.26c}.
Here, we check that each of its hypotheses $(H1),\ldots,(H9)$ is satisfied.
\begin{enumerate}
\item[(H1)] {\it $\bE^1(\Rm)$ is a Fr\'echet space and $\la \bq_k\ra_k$ is a non-decreasing sequence of seminorms that generates its topology.}
\end{enumerate}
\par
{\it Proof.}
See \ref{6.4}.\cqfd
\begin{enumerate}
\item[(H2)] {\it $X[\calG]$ is a locally convex topological vector space whose topology is generated by a filtering family of seminorms $\la \|\cdot\|_i \ra_{i \in I}$ and $\la X_k \ra_k$ is a non-decreasing sequence of $\calG$-closed linear subspaces whose union is $X$.}
\end{enumerate}
\par 
{\it Proof.}
The topology $\calG$ is defined by a single norm $\bG$.
Assume that $\la T_j \ra_j$ is a sequence in $X_k$ that $\bG$-converges to $T \in X$.
The weak* convergence of this sequence to $T$ clearly implies that $\rmspt T \subset B(0,k)$, \ie $T \in X_k$.
Thus $X_k$ is $\bG$-closed.\cqfd
\begin{enumerate}
\item[(H3)] {\it $\lno \cdot \rno$ is a norm on $X$ and $\calG$-lower semicontinuous.}
\end{enumerate}
\par 
{\it Proof.}
See \ref{4.2}(A) and corollary \ref{4.6}.\cqfd
\begin{enumerate}
\item[(H4)] {\it Each $C_k = X_k \cap \{ T \mid \lno T \rno \leq k \}$ is $\calG$-compact.}
\end{enumerate}
\par 
{\it Proof.}
This is theorem \ref{6.2}.\cqfd
\begin{enumerate}
\item[(H5)] {\it $\bD : \bE^1(\Rm) \to X^\star$ is linear.}
\end{enumerate}
\par 
{\it Proof.}
Recalling the definition of $\bD$ in \ref{6.6}, this immediately follows from the linearity of the map $\rmdiv$ in \eqref{eq.4}.\cqfd
\begin{enumerate}
\item[(H6)] {\it For every $v \in \bE^1(\Rm)$ and every $\veps > 0$ there exists $\theta > 0$ such that for every $T \in X_{\lceil \veps^{-1} \rceil}$ we have $|\la T , \bD(v) \ra| \leq \theta \cdot \bG(T) + \veps \cdot \lno T \rno$.}
\end{enumerate}
\par 
{\it Proof.}
Fix $v$.
Given $\veps > 0$, let $k = \lceil \veps^{-1} \rceil$.
Let now $\theta > 0$ be associated with $\veps > 0$ and $k+1$ in theorem \ref{3.3}.
Given $T \in X_k$, choose a sequence $\la u_j \ra_j$ which is $1$-admissible for $T$ at level $k$ (see \ref{6.6}), so that $|\la u_j , \rmdiv v \ra| \leq \theta \cdot \bG(u_j) + \veps \cdot \lno u_j \rno$, by theorem \ref{3.3}.
Therefore, $| \la T , \bD(v) \ra| = \lim_j |\la u_j , \rmdiv v \ra| \leq \theta \cdot \bG(T) + \veps \cdot \lno T \rno$.\cqfd
\begin{enumerate}
\item[(H7)] {\it For every $k \in \N_+$ there exist $\bc_{(H7)}(k) > 0$ and $\kappa(k) \in \N_+$ such that $|\la T , \bD(v) \ra| \leq \bc_{(H7)}(k) \cdot \lno T \rno \cdot \bq_{\kappa(k)}(v)$ whenever $T \in X_k$ and $v \in \bE^1(\Rm)$.}
\end{enumerate}
\par 
{\it Proof.}
We claim this holds with $\kappa(k) = k+2$ and $\bc_{(H7)}(k) = \bc_{\ref{3.2}}(k+1)$.
Indeed, given $v \in \bE^1(\Rm)$ we let $\hat{v}$ be associated with $v$ and $k+1$ in lemma \ref{3.2} and we note that for all $u \in \calD_\#(\Rm)$ with $\rmspt u \subset B(0,k+1)$ we have $|\la u , \rmdiv v \ra| = |\la u , \rmdiv \hat{v} \ra| \leq \lno u \rno \cdot \max_i \|\nabla \hat{v}_i\|_\infty \leq \lno u \rno \cdot \bc_{\ref{3.2}}(k+1) \cdot \bq_{k+2}(v)$.
Thus, if $T \in X_k$ and $\la u_j \ra_j$ is $1$-admissible for $T$ at level $k$ then $|\la T , \bD(v) \ra| = \lim_j |\la u_j , \rmdiv v \ra| \leq \lim_j \lno u_j \rno \cdot \bc_{\ref{3.2}}(k+1) \cdot \bq_{k+2}(v) \leq \lno T \rno \cdot \bc_{\ref{3.2}}(k+1) \cdot \bq_{k+2}(v)$.\cqfd
\begin{enumerate}
\item[(H8)] {\it For every $k \in \N_+$ there exists $\bc_{(H8)}(k) > 0$ such that 
\begin{equation*}
\lno T \rno \leq \bc_{(H8)}(k) \cdot \sup \big\{ |\la T , \bD(v) \ra| \mid v \in \bE^1(\Rm) \text{ and } \bq_k(v) \leq 1 \big\}
\end{equation*}
whenever $T \in X_k$.
}
\end{enumerate}
\par 
{\it Proof.}
We claim this holds with $\bc_{(H8)}(k)=1$.
Let $v \in \calD(\Rm;\Rm)$ be such that $\max_i \|\nabla v_i\|_\infty \leq 1$ and put $w = v - v(0)$.
Then $w \in \bE^1(\Rm)$, $\bq_k(w) = \bq_k(v) \leq 1$, and $\rmdiv w = \rmdiv v \in \calD(\Rm)$.
If $\la u_j \ra_j$ is admissible for $T$ at level $k$ then $\iota(u_j) \to T$ weakly*, since $\bG$-convergence implies weak* convergence, hence
\begin{equation*}
\la T , \bD(w) \ra = \lim_j \int_{\Rm} u_j \cdot \rmdiv v \, d\ssfL^m = \la \rmdiv v , T \ra .
\end{equation*}
Taking the supremum over such $v$, the claim follows from \ref{4.2}(A).\cqfd
\begin{enumerate}
\item[(H9)] {\it For every $k \in \N_+$ and $T \in X \setminus X_k$ there exists $v \in \bE^1(\Rm)$ such that $\la T , \bD(v) \ra \neq 0$ and $\bq_k(v)=0$.}
\end{enumerate}
\par 
{\it Proof.}
Let $k \in \N_+$ and $T \in X \setminus X_k$.
We claim that there exist $i \in \{1,\ldots,m\}$ and $\vphi \in \calD(\Rm)$ such that $\rmspt \vphi \cap B(0,k) = \emptyset$ and $\la D_i\vphi , T \ra \neq 0$.
Assume, if possible, the contrary, \ie $D_iT = 0$ in $\Rm \setminus B(0,k)$ for all $i \in \{1,\ldots,m\}$.
Then $T$ coincides in $\Rm \setminus B(0,k)$ with a constant, by the constancy theorem.
As $\rmspt T$ is compact, this constant is $0$, \ie $\rmspt T \subset B(0,k)$, contradicting $T \not\in X_k$.
Let $i$ and $\vphi$ be as in the claim and put $v = \vphi \cdot e_i$.
Since $0 \not\in \rmspt \vphi$ we have $v \in \calD(\Rm;\Rm) \subset \bE^1(\Rm)$, and $\bq_k(v) = 0$ because $v$ vanishes in a neighbourhood of $B(0,k)$.
Finally, $\la T , \bD(v) \ra = \la \rmdiv v , T \ra = \la D_i\vphi , T \ra \neq 0$, as in the proof of (H8).\cqfd
\end{Empty}
\begin{Theorem}
\label{6.8}
Let $F : \calD_\#(\Rm) \to \R$ be linear.
The following are equivalent.
\begin{enumerate}
\item[(A)] There exists $v \in C^1(\Rm;\Rm)$ such that $\la u , F \ra = \int_{\Rm} u \cdot \rmdiv v \, d\ssfL^m$ for all $u \in \calD_\#(\Rm)$.
\item[(B)] $(\forall \veps > 0)(\forall k \in \N)(\exists \theta > 0)(\forall u \in \calD_\#(\Rm)) : \rmspt u \subset B(0,k) \, \Rightarrow \, |\la u , F \ra | \leq \theta \cdot \bG(u) + \veps \cdot \lno u \rno$.
\item[(C)] There exists a $\calG_{\lno \cdot \rno}$-continuous linear functional $\hat{F} : X \to \R$ such that $\hat{F} \circ \iota = F$.
\end{enumerate}
Furthermore, under these conditions, if $\veps > 0$ and $k \in \N_+$ is so that $\la u , F \ra \neq 0$ for some $u \in \calD_\#(\Rm)$ with $\rmspt u \subset B(0,k)$ then $v$ satisfying (A) can be chosen such that also
\begin{equation*}
\bq_k(v) \leq (1 + \veps) \cdot \vvvert F \vvvert_{k+1},
\end{equation*}
where
\begin{equation*}
\vvvert F \vvvert_k = \sup \left\{ |\la u , F \ra| \mid u \in \calD_\#(\Rm), \ \rmspt u \subset B(0,k), \text{ and } \lno u \rno \leq 1 \right\}.
\end{equation*}
\end{Theorem}

\begin{proof}
{(A) $\Rightarrow$ (B).}
This is theorem \ref{3.3}.
\par
{(B) $\Rightarrow$ (C).}
Regard $Y=\iota(\calD_\#(\Rm))$ as a linear subspace of $X$, according to \ref{4.2}(E), and consider $F_Y : Y \to \R$ defined by $F_Y(\iota(u))=\la u,F\ra$.
Given $\veps>0$, condition (B), applied with $k=\lceil\veps^{-1}\rceil$, shows that $F_Y$ satisfies hypothesis (H7) of \cite[Theorem 5.7]{DEP.26c}.
The other hypotheses of that theorem follow from \ref{6.1}, \ref{6.7}(H2),(H3), and the uniform sequential density established in \ref{6.5}.
Thus $F_Y$ admits a $\calG_{\lno\cdot\rno}$-continuous linear extension $\hat{F}:X\to\R$, and $\hat{F}\circ\iota=F$.
\par
{(C) $\Rightarrow$ (A).}
By \ref{6.7}, the hypotheses of \cite[Theorem 8.1]{DEP.26c} are satisfied, so its conclusion (C2) yields $v\in\bE^1(\Rm)$ such that $\bD(v)=\hat{F}$.
For every $u\in\calD_\#(\Rm)$, the definition of $\bD$ in \ref{6.6} gives
\begin{equation*}
\la u,F\ra=\la\iota(u),\hat{F}\ra=\la\iota(u),\bD(v)\ra=\int_{\Rm}u\cdot\rmdiv v\,d\ssfL^m.
\end{equation*}
Since $v\in\bE^1(\Rm)\subset C^1(\Rm;\Rm)$, this proves (A).
\par 
Finally, assume that $\veps > 0$ and that $k \in \N_+$ satisfies the hypothesis specified after ``Furthermore''. 
By \ref{6.5}, we have
\begin{equation*}
\sup \left\{ |\la T , \hat{F} \ra| \mid T \in X_k \text{ and } \lno T \rno \leq 1\right\} \leq \vvvert F \vvvert_{k+1}.
\end{equation*}
The supremum on the left is positive by hypothesis.
The conclusion follows from \cite[Theorem 8.1(C3)]{DEP.26c}, as the constant $\bc_{(H8)}(k)=1$ in \ref{6.7}(H8).
\end{proof}

\begin{Empty}[Awkwardness of $\calG_{\lno \cdot \rno}$]
\label{6.9}
Recall that the topology $\calG_{\lno \cdot \rno}$ of $X$ is defined in \ref{6.1} and observe its relevance: theorem \ref{6.8} identifies the divergence of a $C^1$ vector field with a member of the dual space $X[\calG_{\lno \cdot \rno}]^*$.
\begin{enumerate}
\item[(A)] {\it The topology $\calG_{\lno \cdot \rno}$ is sequential.}
\end{enumerate}
\par
{\it Proof.}
Each $C_k$ is $\bG$-compact, by theorem \ref{6.2}, and the topology induced by $\bG$ on $C_k$ is metrisable.
The assertion follows from \cite[Theorem 4.4(C1)]{DEP.26c}.\cqfd
\par
We next show that $\lno\cdot\rno$ is not $\calG_{\lno\cdot\rno}$-continuous.
Let $e_1$ be the first vector of the canonical basis of $\Rm$, put $T_0=\bdelta_{e_1}-\bdelta_0$, and abbreviate $s=\lno T_0\rno>0$.
For $0<r\leq1$, define $T_r=s^{-1}\cdot\bmu_{r\#}T_0=s^{-1}\cdot(\bdelta_{r\cdot e_1}-\bdelta_0)$.
According to \ref{4.3}(A) and (B),
\begin{equation*}
\bG(T_r)=\frac{r}{s}\cdot\bG(T_0) \qquad\text{and}\qquad \lno T_r\rno=1.
\end{equation*}
As $T_r\in C_1$, it follows from \cite[Theorem 3.1(C1)]{DEP.26c} that $T_r\to0$ with respect to $\calG_{\lno\cdot\rno}$ as $r\to0$, whereas $\lno T_r\rno=1$.
\begin{enumerate}
\item[(B)] {\it The topology $\calG_{\lno \cdot \rno}$ is not barrelled.}
\end{enumerate}
\par
{\it Proof.}
By \cite[Theorem 3.4(C3)]{DEP.26c}, the norm $\lno\cdot\rno$ is $\calG_{\lno\cdot\rno}$-lower-semicontinuous.
If $X[\calG_{\lno\cdot\rno}]$ were barrelled, this norm would be continuous, see \eg \cite[D.4(b)]{DEP.26c}, contradicting the preceding observation.\cqfd
\begin{enumerate}
\item[(C)] {\it The topology $\calG_{\lno \cdot \rno}$ is not bornological.}
\end{enumerate}
\par
{\it Proof.}
By \cite[Theorem 3.4(C3)]{DEP.26c}, the norm $\lno\cdot\rno$ is bounded on each $\calG_{\lno\cdot\rno}$-bounded set.
If $X[\calG_{\lno\cdot\rno}]$ were bornological, this norm would be continuous, see \eg \cite[E.2]{DEP.26c}, a contradiction.\cqfd
\begin{enumerate}
\item[(D)] {\it The topology $\calG_{\lno \cdot \rno}$ is not Fr\'echet-Urysohn.}
\end{enumerate}
\par
{\it Proof.}
This follows from (C), since every Fr\'echet-Urysohn locally convex topological vector space is bornological, see \eg \cite[E.5]{DEP.26c}.\cqfd
\end{Empty}

\section{Necessary and sufficient condition for membership in $\rmdiv(\rmLip(\Rm;\Rm))$}

\begin{Empty}[The norm $ \bp$]
\label{7.1}
We set
\begin{equation*}
\rmLip_{0,1}(\Rm;\Rm) = \left\{ v : \Rm \to \Rm \mid v(0)=0 \text{ and } \max_i \|\nabla v_i\|_{L_\infty} \leq 1\right\}.
\end{equation*}
For $f \in L_1(\Rm)$ we put
\begin{equation*}
\bp(f) = \sup \left\{ \int_{\Rm} f \cdot \rmdiv v \, d\ssfL^m \mid v \in \rmLip_{0,1}(\Rm;\Rm) \right\}.
\end{equation*}
\begin{enumerate}
\item[(A)] {\it $\bp$ is a norm on $L_1(\Rm)$ and $\bp \leq m \cdot \|\cdot\|_{L_1}$.}
\end{enumerate}
\par 
{\it Proof.}
Let $f \in L_1(\Rm)$.
If $v \in \rmLip_{0,1}(\Rm;\Rm)$ then $\|\rmdiv v\|_{L_\infty} \leq m$, thus, $f \cdot \rmdiv v \in L_1(\Rm)$ so that $\bp(f)$ is well-defined and $\bp(f) \leq m \cdot \|f\|_{L_1}$.
It is now clear that $\bp$ is a semi-norm.
If $\bp(f)=0$ then using $v = \vphi \cdot e_i$ for all $\vphi \in \calD(\Rm)$ we infer that $D_if=0$ distributionally for all $i \in \{1,\ldots,m\}$.
By the constancy theorem and the integrability of $f$, we infer that $f=0$.\cqfd
\begin{enumerate}
\item[(B)] {\it $\left| \int_{\Rm} f \, d\ssfL^m \right| \leq \frac{1}{m} \cdot \bp(f)$ for all $f \in L_1(\Rm)$.}
\end{enumerate}
\par 
{\it Proof.}
Use $v = \pm \rmid_{\Rm}$ in the definition of $\bp$.\cqfd
\begin{enumerate}
\item[(C)] {\it $\lno u \rno \leq \bp(u) \leq 5 \cdot \lno u \rno$ for all $u \in \calD_\#(\Rm)$.}
\end{enumerate}
\par 
{\it Proof.}
The first inequality is a consequence of \ref{3.1}(A) and the fact that $v - v(0) \in \rmLip_{0,1}(\Rm;\Rm)$ whenever $v \in \calD(\Rm;\Rm)$ and $\max_i \|\nabla v_i\|_\infty \leq 1$.
Conversely, let $u \in \calD_\#(\Rm)$ and $r > 0$ be such that $\rmspt u \subset U(0,r)$.
Choose $\psi \in \calD(\Rm)$ such that $\ind_{B(0,r)} \leq \psi \leq \ind_{B(0,2r)}$ and $\|\nabla \psi\|_\infty \leq \frac{2}{r}$. Given $v \in \rmLip_{0,1}(\Rm;\Rm)$ we define a compactly supported Lipschitz vector field $w = \psi \cdot v$.
For each $i \in \{1,\ldots,m\}$ we have $\nabla w_i = v_i \cdot \nabla \psi + \psi \cdot \nabla v_i$, hence, $\|\nabla w_i\|_{L_\infty} \leq 2r \cdot \frac{2}{r} + 1 = 5$, as $|v_i(x)| = |v_i(x)-v_i(0)| \leq 2r$ whenever $x \in \rmspt \nabla \psi$.
Since $w = v$ on a neighbourhood of $\rmspt u$ we have $\int_{\Rm} u \cdot \rmdiv v \, d\ssfL^m = \int_{\Rm} u \cdot \rmdiv w \, d\ssfL^m$.
For each $\delta > 0$ we have $w * \Phi_\delta \in \calD(\Rm;\Rm)$ and $\max_i \|\nabla [w * \Phi_\delta]_i\|_\infty \leq 5$, therefore, $\int_{\Rm} u \cdot \rmdiv (w * \Phi_\delta) \, d\ssfL^m \leq 5 \cdot \lno u \rno$.
Since $\lim_{\delta \to 0^+} \int_{\Rm} u \cdot \rmdiv (w * \Phi_\delta) \, d\ssfL^m = \int_{\Rm} u \cdot \rmdiv v \,d\ssfL^m$, the conclusion follows from the arbitrariness of $v$.\cqfd
\begin{enumerate}
\item[(D)] {\it $\bp(f \circ \btau_a) = \bp(f) = \bp(r^m \cdot f \circ \bmu_r)$ for all $f \in L_1(\Rm)$, $a \in \Rm$, and $r > 0$.}
\end{enumerate}
\par 
{\it Proof.}
Letting $\hat{f}(x) = (f \circ \btau_a)(x) = f(x+a)$ and $\hat{v} \in \rmLip_{0,1}(\Rm;\Rm)$ we note that $\int_{\Rm} \hat{f} \cdot \rmdiv \hat{v} \, d\ssfL^m = \int_{\Rm} f \cdot \rmdiv v \, d\ssfL^m$, where $v(x) = \hat{v}(x-a) - \hat{v}(-a)$ so that $v \in \rmLip_{0,1}(\Rm;\Rm)$.
Similarly, letting $\hat{f}(x) = r^m \cdot (f \circ \bmu_r)(x) = r^m \cdot f(r \cdot x)$ and $\hat{v} \in \rmLip_{0,1}(\Rm;\Rm)$ we observe that $\int_{\Rm} \hat{f} \cdot \rmdiv \hat{v} \,d\ssfL^m = \int_{\Rm} f \cdot \rmdiv v \, d\ssfL^m$, where $v(x) = r \cdot \hat{v} \left( \frac{x}{r} \right)$ so that $v \in \rmLip_{0,1}(\Rm;\Rm)$.\cqfd
\par 
Finally, we define
\begin{equation*}
K = L_\infty(\Rm) \cap \left\{ \rmdiv v \mid v \in \rmLip_{0,1}(\Rm;\Rm) \right\}.
\end{equation*}
\end{Empty}

\begin{Theorem}
\label{7.2}
The set $K$ is absolutely convex, weakly*-compact, and
\begin{equation}
\label{eq.7.1}
K = L_\infty(\Rm) \cap \left\{ g \mid \int_{\Rm} f \cdot g \, d\ssfL^m \leq \bp(f) \text{ for all } f \in L_1(\Rm) \right\}.
\end{equation}
\end{Theorem}

The above is a necessary and sufficient condition for a function $g \in L_\infty(\Rm)$ to be the divergence of a Lipschitz vector field in $\Rm$.
We give three (similar) proofs below.

\begin{proof}[First proof]
{\bf (i)}
The absolute convexity of $K$ is obvious.
As $K$ is contained in $B_{L_\infty}(0,m)$, its weak* compactness will follow, by Banach-Alaoglu's theorem, from it being weakly*-closed.
Since $L_1(\Rm)$ is separable, the restriction to $B_{L_\infty}(0,m)$ of the weak* topology is metrisable.
Hence, it suffices to show that each sequence $\la g_j \ra_j$ in $K$ admits a subsequence converging weakly* to a member of $K$.
Let $v_j \in \rmLip_{0,1}(\Rm;\Rm)$ be so that $g_j = \rmdiv v_j$.
By Ascoli's theorem and a diagonal argument on an exhaustion by compact sets, some subsequence -- still denoted $\la v_j \ra_j$ -- converges locally uniformly to some $v \in \rmLip_{0,1}(\Rm;\Rm)$.
If $\vphi \in \calD(\Rm)$ then $\lim_j \int_{\Rm} \vphi \cdot g_j \, d\ssfL^m = - \lim_j \int_{\Rm} (\nabla \vphi) \ip v_j \,d\ssfL^m = - \int_{\Rm} (\nabla \vphi) \ip v \,d\ssfL^m = \int_{\Rm} \vphi \cdot g \, d\ssfL^m$, where $g = \rmdiv v \in K$.
Since $\calD(\Rm)$ is dense in $L_1(\Rm)$ and $\la g_j \ra_j$ is bounded in $L_\infty(\Rm)$, we easily conclude that $\la g_j \ra_j$ converges weakly* to $g$.
\par 
{\bf (ii)}
We turn to proving \eqref{eq.7.1}.
The inclusion $\subset$ immediately follows from the definitions of $K$ and $\bp$.
We prove the reverse inclusion by contraposition.
If $g_0 \not \in K$ then Hahn-Banach's theorem \cite[3.7]{RUDIN} applies to $K$ in $L_\infty(\Rm)[\sigma(L_\infty,L_1)]$: There exists $f \in L_1(\Rm)$ such that $\int_{\Rm} f \cdot g_0 \, d\ssfL^m > \sup \left\{ \int_{\Rm} f\cdot g \, d\ssfL^m \mid g \in K \right\} = \bp(f)$.
\end{proof}

\begin{proof}[Second proof]
We retain part {\bf (i)} of the proof above and simply note that part {\bf (ii)} is a proof of the bipolar theorem applied to $K$, see \eg \cite[C.4(C)]{DEP.26c}.
\end{proof}

\begin{proof}[Third proof]
We dispense with part {\bf (i)} of the first proof by using duality with the Arens-Eells (or Lipschitz-free) space $\calF(\Rm)$ (we refer to \cite[§3.1]{DEP.26a} for notation and facts).
First note that if $f \in L_1(\Rm)$ and $i \in \{1,\ldots,m\}$ then $D_if$ may be interpreted as a member of $\calF(\Rm)$, according to \cite[3.1.3(B)]{DEP.26a} applied to $v = f \cdot e_i$, \ie the function $\rmLip_0(\Rm) \to \R : u \mapsto - \int_{\Rm} f \cdot D_i u \, d\ssfL^m$ is a member of the predual $\calF(\Rm)$ of $\rmLip_0(\Rm)$ and$\|D_if\|_\calF \leq \|f\|_{L_1}$.
Therefore, the following linear map is bounded 
\begin{equation*}
\nabla : L_1(\Rm) \to \bigoplus_{\ell_1}^m \calF(\Rm) : f \mapsto D_1f \oplus \cdots \oplus D_m f.
\end{equation*}
In view of the next to last sentence, its adjoint is
\begin{equation*}
-\rmdiv : \bigoplus_{\ell_\infty}^m \rmLip_0(\Rm) \to L_\infty(\Rm) : v_1 \oplus \cdots \oplus v_m \mapsto -\sum_{i=1}^m D_iv_i.
\end{equation*}
Since $B$, the unit ball of the domain, is symmetric, $K = \rmdiv(B) = (-\rmdiv)(B)$.
The weak*-compactness of $K$ follows from that of $B$ (Banach-Alaoglu) and the weak*-to-weak* continuity of $-\rmdiv$.
\end{proof}

\begin{Empty}[$\bp_\veps$ and $\bp^*_\veps$]
\label{7.3}
For $\veps > 0$ and $f \in L_1(\Rm)$ we define
\begin{equation*}
\bp_\veps(f) = \sup \left\{ \int_{\Rm} f \cdot \rmdiv v \, d\ssfL^m \mid v \in \rmLip_{0,1}(\Rm;\Rm) \text{ and } \|\rmdiv v\|_{L_\infty} \leq \veps \right\}
\end{equation*}
and 
\begin{equation*}
\bp^*_\veps(f) = \inf \left\{ \veps \cdot \|h\|_{L_1} + \bp(f-h) \mid h \in L_1(\Rm) \right\}.
\end{equation*}
\begin{enumerate}
\item[(A)] {\it $\bp^*_\veps$ is a seminorm on $L_1(\Rm)$.}
\end{enumerate}
\par 
{\it Proof.}
Trivial.\cqfd
\begin{enumerate}
\item[(B)] {\it $\bp^*_\veps(f) \leq \veps \cdot \|f\|_{L_1}$ and $\bp^*_\veps(f) \leq \bp(f)$ for all $f \in L_1(\Rm)$.}
\end{enumerate}
\par
{\it Proof.}
Take, respectively, $h=f$ and $h=0$.\cqfd
\end{Empty}

\begin{Proposition}
\label{7.4}
$\bp_\veps = \bp^*_\veps$ and the supremum in the definition of $\bp_\veps$ is attained.
\end{Proposition}

\begin{proof}
{\bf (i)}
We first show that $\bp_\veps \leq \bp^*_\veps$.
Let $f \in L_1(\Rm)$ and $v \in \rmLip_{0,1}(\Rm;\Rm)$ be such that $\|\rmdiv v \|_{L_\infty} \leq \veps$.
For each $h \in L_1(\Rm)$ we have $\int_{\Rm} f \cdot \rmdiv v \,d\ssfL^m = \int_{\Rm} h \cdot \rmdiv v \,d\ssfL^m + \int_{\Rm} (f-h) \cdot \rmdiv v \, d\ssfL^m \leq \veps \cdot \|h\|_{L_1} + \bp(f-h)$.
As $f$ and $h$ are arbitrary, $\bp_\veps(f) \leq \bp^*_\veps(f)$.
\par 
{\bf (ii)}
Let $f_0 \in L_1(\Rm)$.
By \ref{7.3}(A) and Hahn-Banach's theorem \cite[3.2]{RUDIN}, there exists a linear $\alpha : L_1(\Rm) \to \R$ such that $\alpha(f_0) = \bp^*_\veps(f_0)$ and $|\alpha| \leq \bp^*_\veps$.
As $\bp^*_\veps \leq \veps \cdot \|\cdot\|_{L_1}$,  by \ref{7.3}(B), there exists $g \in L_\infty(\Rm)$ such that $\alpha(f) = \int_{\Rm} f \cdot g \, d\ssfL^m$ for every $f \in L_1(\Rm)$.
By \ref{7.3}(B) again, $\|g\|_{L_\infty} \leq \veps$ and $\int_{\Rm} f \cdot g \,d\ssfL^m = \alpha(f) \leq \bp(f)$ for every $f \in L_1(\Rm)$.
Thus, $g = \rmdiv v$ for some $v \in \rmLip_{0,1}(\Rm;\Rm)$, according to theorem \ref{7.2}.
Finally, since $\|\rmdiv v\|_{L_\infty} = \|g\|_{L_\infty} \leq \veps$, the vector field $v$ is admissible in the definition of $\bp_\veps$, and hence $\bp^*_\veps(f_0) = \alpha(f_0) = \int_{\Rm} f_0 \cdot \rmdiv v \, d\ssfL^m \leq \bp_\veps(f_0)$. Combined with {\bf (i)}, this proves both the equality and attainment of the supremum.
\end{proof}

\section{$\rmdiv(C^1(\Rm;\Rm)) \subsetneqq \rmdiv(\rmLip_{\rmloc}(\Rm;\Rm)) \cap C(\Rm)$}

\begin{Theorem}[Dual to Ornstein's non-inequality]
\label{8.1}
For every $\veps > 0$ there exists $\vphi \in \calD(\Rm)$ with $\bp(\vphi) \leq 1$ and there exists $v \in \rmLip_{0,1}(\Rm;\Rm)$ with $\|\rmdiv v \|_{L_\infty} \leq \veps$ such that $\int_{\Rm} \vphi \cdot \rmdiv v \, d\ssfL^m \geq \frac{1}{2}$.
\end{Theorem}

\begin{proof}
{\bf (i)} 
The proof is by contradiction, \ie we assume that there exists $\veps > 0$ such that for all $\vphi \in \calD(\Rm)$ with $\bp(\vphi) \leq 1$ and all $v \in \rmLip_{0,1}(\Rm;\Rm)$ with $\|\rmdiv v\|_{L_\infty} \leq \veps$ we have $\int_{\Rm} \vphi \cdot \rmdiv v \, d\ssfL^m < \frac{1}{2}$.
Recalling the definition of $\bp_\veps$ and the homogeneity of $\bp$, this implies that $\bp_\veps(\vphi) \leq \frac{1}{2} \cdot \bp(\vphi)$ for all $\vphi \in \calD(\Rm)$ and, since $\calD(\Rm)$ is dense in $L_1(\Rm)$, that $\bp_\veps(f) \leq \frac{1}{2} \cdot \bp(f)$ for all $f \in L_1(\Rm)$.
It then ensues from proposition \ref{7.4} that for every $f \in L_1(\Rm)$ there exists $h \in L_1(\Rm)$ such that
\begin{equation}
\label{eq.8.1}
\veps \cdot \|h\|_{L_1} + \bp(f-h) \leq \frac{3}{4} \cdot \bp(f).
\end{equation}
\par 
{\bf (ii)}
We claim that $\|f\|_{L_1} \leq \frac{3}{\veps} \cdot \bp(f)$ for all $f \in L_1(\Rm)$.
Given $f$ we define by induction a sequence $\la f_j \ra_j$ in $L_1(\Rm)$ by letting $f_1=f$ and $f_{j+1} = f_j - h_j$, where $h_j$ is associated with $f_j$ in \eqref{eq.8.1}.
Observe that $\bp(f_{j+1}) \leq \frac{3}{4} \cdot \bp(f_j) \leq \left( \frac{3}{4}\right)^j \cdot \bp(f)$ and $\|h_{j+1}\|_{L_1} \leq \frac{1}{\veps} \cdot \frac{3}{4} \cdot \bp(f_j) \leq \frac{1}{\veps} \cdot \left( \frac{3}{4} \right)^j \cdot \bp(f)$ for all $j$. 
In particular, the series $\sum_{j=1}^\infty h_j$ converges absolutely in $L_1(\Rm)$ to some $h$. 
Letting $n \to \infty$ in the identity $f_{n+1} = f - \sum_{j=1}^n h_j$ we infer that the sequence $\la f_j \ra_j$ converges to $f-h$ in $L_1(\Rm)$.
In view of \ref{7.1}(A), it ensues that $\bp(f-h) = \lim_j \bp(f_j) = 0$.
As $\bp$ is a norm, $f=h$.
Finally, $\|f\|_{L_1} = \|h\|_{L_1} \leq \sum_{j=1}^\infty \|h_j\|_{L_1} \leq \frac{1}{\veps} \cdot \left( \sum_{j=1}^\infty \left( \frac{3}{4}\right)^j\right) \cdot \bp(f)$.
\par 
{\bf (iii)}
It remains to note that {\bf (ii)} contradicts Ornstein's non-inequality.
Indeed, letting $\bc = \frac{15 \cdot m}{\veps\cdot (m-1)}$ and applying \ref{Ornstein}(B) we find $\vphi \in \calD(\Rm)$ such that $\|D^2_0 \vphi\|_{L_1}=1$ and $\bc < \|\triangle \vphi\|_{L_1}$. 
Putting $u = \triangle \vphi$ we see that $u \in L_1(\Rm)$ and, by {\bf (ii)}, \ref{7.1}(C), and Theorem \ref{4.2F} we have $\bc < \|u\|_{L_1} \leq \frac{3}{\veps} \cdot \bp(u) \leq \frac{15}{\veps} \cdot \lno u \rno \leq \frac{15}{\veps} \cdot \frac{m}{m-1} = \bc$.
\end{proof}

\begin{Remark}
\label{8.2}
Let $\vphi$ and $w$ be as in theorem \ref{8.1} and choose $a \in \Rm$ and $r > 0$ such that $\rmspt \vphi \subset U(a,r)$.
Define $\hat{\vphi}(x) = r^m \cdot \vphi \left( a + r \cdot x \right)= r^m \cdot (\vphi \circ \btau_a \circ \bmu_r)(x)$ and $\hat{w}(x) = \frac{1}{r} \cdot w \left( a+r \cdot x\right) - \frac{1}{r} \cdot w \left( a \right)$, $x \in \Rm$.
Note that $\rmspt \hat{\vphi} \subset U(0,1)$, $\hat{\vphi} \in \calD(\Rm)$, and $\bp(\hat{\vphi}) = \bp(\vphi) \leq 1$, by \ref{7.1}(D).  
Note also that $\hat{w}(0)=0$ and $\rmLip \hat{w} = \rmLip w$ so that $\hat{w} \in \rmLip_{0,1}(\Rm)$.
Moreover, $(\rmdiv \hat{w})(x) = (\rmdiv w) \left(a + r \cdot x \right)$ for all $x \in \Rm$ so that $\|\rmdiv \hat{w}\|_{L_\infty} = \|\rmdiv w\|_{L_\infty} \leq \veps$.
Finally, $\int_{\Rm} \hat{\vphi} \cdot \rmdiv \hat{w} \, d\ssfL^m = \int_{\Rm} \vphi \cdot \rmdiv w \, d\ssfL^m \geq \frac{1}{2}$.
In other words, we may add as a conclusion of theorem \ref{8.1} that $\rmspt \vphi \subset U(0,1)$.
\end{Remark}

\begin{Corollary}[Building block]
\label{8.3}
There exists $\bc_{\theTheorem}(m) > 0$ with the following property.
For every $0 < \veps \leq \frac{1}{8}$ there are $u \in \calD_\#(\Rm)$ and $v \in C^\infty(\Rm;\Rm)$ such that:
\begin{enumerate}
\item[(A)] $\rmspt u \subset U(0,1)$ and $\lno u \rno \leq 2$;
\item[(B)] $\rmspt (\rmdiv v) \subset U(0,2)$;
\item[(C)] $\|\rmdiv v\|_{L_\infty} \leq \veps$;
\item[(D)] $\displaystyle \int_{\Rm} u \cdot \rmdiv v \, d\ssfL^m \geq \frac{3}{8}$;
\item[(E)] $\|v\|_\infty \leq \bc_{\theTheorem}(m)$ and $\max_i \|\nabla v_i\|_\infty \leq \bc_{\theTheorem}(m)$;
\item[(F)] For all $x \in \Rm \setminus B(0,4)$ and all $i \in \{1,\ldots,m\}$ we have $\displaystyle|\nabla v_i(x)|_2 \leq \frac{\bc_{\theTheorem}(m)}{|x|_2^m}$.
\end{enumerate}
\end{Corollary}

\begin{proof}
{\bf (i)} 
Fix $0 < \veps \leq \frac{1}{8}$.
Let $\vphi$ and $w$ be given by theorem \ref{8.1} so that $\bp(\vphi) \leq 1$, $\|\rmdiv w \|_{L_\infty} \leq \veps$, $\int_{\Rm} \vphi \cdot \rmdiv w \, d\ssfL^m \geq \frac{1}{2}$, and $\rmspt \vphi \subset U(0,1)$, the latter being a consequence of remark \ref{8.2}.
Fix $\theta \in \calD(\Rm)$ such that $\rmspt \theta \subset U(0,1)$, $\theta \geq 0$, and $\int_{\Rm} \theta \,d\ssfL^m = 1$.
Put $a = \int_{\Rm} \vphi \,d\ssfL^m$ so that $|a| \leq \frac{1}{m}$, by \ref{7.1}(B), and define $u = \vphi - a \cdot \theta$.
Clearly, $u \in \calD_\#(\Rm)$ and $\rmspt u \subset U(0,1)$.
Moreover, $\lno u \rno \leq \bp(u) \leq \bp(\vphi) + |a| \cdot m \cdot \|\theta\|_{L_1} \leq 2$, by \ref{7.1}(C,A).
This completes the definition of $u$ and the proof of conclusion (A).
Finally, note that
\begin{equation}
\label{eq.8.3.pairing}
\int_{\Rm} u \cdot \rmdiv w \, d\ssfL^m = \int_{\Rm} \vphi \cdot \rmdiv w \, d\ssfL^m - a \cdot \int_{\Rm} \theta \cdot \rmdiv w \, d\ssfL^m \geq \frac{7}{16}
\end{equation}
since $|a| \leq \frac{1}{2}$ and $\left| \int_{\Rm} \theta \cdot \rmdiv w \, d\ssfL^m \right| \leq \|\theta\|_{L_1} \cdot \|\rmdiv w \|_{L_\infty} \leq \frac{1}{8}$.
\par 
{\bf (ii)}
We ought to replace $w$ by a smooth vector field whose divergence is compactly supported in $U(0,2)$ and whose divergence is close enough to that of $w$ near $\rmspt u$ for (D) to hold.
For $\delta > 0$ we abbreviate $w_\delta = w * \Phi_\delta$ so that $w_\delta \in C^\infty(\Rm;\Rm)$, $\|\rmdiv w_\delta\|_\infty = \|(\rmdiv w)*\Phi_\delta\|_\infty \leq \|\rmdiv w \|_{L_\infty} \leq \veps$ and 
\begin{equation}
\label{eq.8.2}
\int_{\Rm} u \cdot \rmdiv w_\delta \, d\ssfL^m \geq \frac{3}{8}
\end{equation}
provided $\delta$ is sufficiently small.
From now on, such $0 < \delta \leq 1$ is fixed.
We observe that 
\begin{equation}
\label{eq.8.3}
\| \nabla (w_\delta \ip e_i) \|_\infty = \| (\nabla w_i) * \Phi_\delta \|_\infty \leq \| \nabla w_i \|_{L_\infty} \leq 1,
\end{equation}
since $w \in \rmLip_{0,1}(\Rm)$.
In particular, $\rmLip(w_\delta) \leq \sqrt{m}$.
As $w(0)=0$ and also $\rmLip w \leq \sqrt{m}$ we infer that $|w_\delta(0)|_2 \leq \sqrt{m} \cdot \delta$ and, in turn,
\begin{equation}
\label{eq.8.4}
|w_\delta(x)|_2 \leq \sqrt{m} \cdot (2 + \delta) \leq 3 \cdot \sqrt{m}
\end{equation}
for all $x \in U(0,2)$.
\par 
Choose a radial $\psi \in \calD(\Rm)$ such that $\ind_{B(0,1)} \leq \psi \leq \ind_{B(0,3/2)}$, $\|\nabla \psi \|_\infty \leq 3$, and $\|D^2\psi\|_{\infty,\mathrm{op}} \leq 24$. Such a cutoff is obtained by smoothing a radial $C^{1,1}$ cutoff with these bounds.
Define $\hat{w} = \psi \cdot w_\delta$ and define smooth functions $f$ and $f_{\err}$ as follows:
\begin{equation}
\label{eq.8.5}
\rmdiv \hat{w} = \underbrace{\psi \cdot \rmdiv w_\delta}_{f} + \underbrace{(\nabla \psi) \ip w_\delta}_{f_{\err}}.
\end{equation}
Regarding $f$, we observe that
\begin{equation}
\label{eq.8.6}
\rmspt f \subset U(0,2) \,, \quad \|f\|_\infty \leq \veps \,,\quad\text{and} \quad \int_{\Rm} u \cdot f \,d\ssfL^m \geq \frac{3}{8}\,,
\end{equation}
the latter being a consequence of \eqref{eq.8.2}, since $\rmspt u \subset \{\psi = 1 \}$.
About $f_{\err}$ we note that 
\begin{equation}
\label{eq.8.7}
\|f_{\err}\|_\infty \leq \|\nabla \psi\|_\infty \cdot \sup_{x \in U(0,2)} |w_\delta(x)|_2 \leq 9 \cdot \sqrt{m},
\end{equation}
by \eqref{eq.8.4}. Furthermore, differentiating $f_{\err} = (\nabla\psi)\ip w_\delta$, we obtain $\nabla f_{\err} = (D^2\psi)\cdot w_\delta + (Dw_\delta)^T\cdot\nabla\psi$. Thus, by \eqref{eq.8.3}, \eqref{eq.8.4}, and the bounds on $\psi$, we have $\|\nabla f_{\err}\|_\infty \leq 24\cdot 3\sqrt{m} + 3\sqrt{m}$, therefore,
\begin{equation}
\label{eq.8.8}
\rmLip f_{\err} \leq 75 \cdot \sqrt{m}.
\end{equation}
\par 
{\bf (iii)}
In order to define $v$, we need to write $f_{\err}$ as a divergence.
Applying \ref{Newton} (with $f=f_{\err}$, recalling that $\rmspt f_{\err} \subset B(0,2)$) define $w_{\err} = (\nabla \Gamma) * f_{\err} \in C^\infty(\Rm;\Rm)$ so that $\rmdiv w_{\err} = f_{\err}$.
Consequently, upon letting $v = \hat{w} - w_{\err}$ we note that $v \in C^\infty(\Rm;\Rm)$, $\rmdiv v = f$, by \eqref{eq.8.5}, so that conclusions (B), (C), and (D) immediately follow from \eqref{eq.8.6}.
\par 
We notice that 
\begin{equation}
\label{eq.8.9}
\|v\|_\infty \leq \|\hat{w}\|_\infty + \|w_{\err}\|_\infty \leq \sup_{x \in U(0,2)} |w_\delta(x)|_2 + 2 \cdot \|f_{\err}\|_\infty \leq 3 \cdot \sqrt{m} + 18 \cdot \sqrt{m},
\end{equation}
by \eqref{eq.8.4}, \ref{Newton}(B), and \eqref{eq.8.7}.
Furthermore, $\|\nabla (\hat{w} \ip e_i)\|_\infty \leq \|\nabla \psi\|_\infty \cdot \sup_{x \in U(0,2)} |(w_\delta \ip e_i)(x)| + \|\psi\|_\infty \cdot \|\nabla (w_\delta \ip e_i)\|_\infty \leq 9 \cdot \sqrt{m} + 1$, by the definition of $\hat{w}$, \eqref{eq.8.4}, and \eqref{eq.8.3}.
We also infer from \ref{Newton}(A)(a) that $\|D_j(w_{\err} \ip e_i)\|_\infty = \|D^2_{ij}(\Gamma * f_{\err})\|_\infty \leq (3m - 2)\cdot 2 \cdot \rmLip f_{\err} \leq 150 \cdot (3m-2) \cdot \sqrt{m}$, by \eqref{eq.8.8}.
Thus,
\begin{equation}
\label{eq.8.10}
\|\nabla(v \ip e_i)\|_\infty \leq \|\nabla (\hat{w} \ip e_i)\|_\infty + \|\nabla (w_{\err} \ip e_i)\|_\infty \leq 9 \cdot \sqrt{m} + 1 + 150 \cdot m \cdot (3m-2) \cdot \sqrt{m}.
\end{equation}
This establishes conclusion (E) provided that $\bc_{\ref{8.3}}(m) \geq 1 + 21 \cdot \sqrt{m} + 150 \cdot m \cdot (3m-2) \cdot \sqrt{m}$.
\par 
It remains to prove conclusion (F).
First notice that $\|f_{\err}\|_{L_1} \leq2^m \cdot \balpha(m) \cdot 9 \cdot \sqrt{m}$, by \eqref{eq.8.7} and $\rmspt f_{\err} \subset B(0,2)$.
If $|x|_2 > 4$ then $D_j(v \ip e_i)(x) = -D_j(w_{\err} \ip e_i)(x) = -[D^2_{ij}(\Gamma*f_{\err})](x)$, since $\rmspt \hat{w} \subset B(0,3/2)$, so that, according to \ref{Newton}(A)(b), we have
\begin{equation*}
|D_j(v \ip e_i)(x)| \leq \frac{(m-1) \cdot 2^m}{m \cdot \balpha(m)} \frac{\|f_{\err}\|_{L_1}}{|x|_2^m} \leq \frac{(m-1) \cdot 4^m}{m} \cdot \frac{9 \cdot \sqrt{m}}{|x|_2^m}.
\end{equation*}
Therefore, (F) holds provided that $\bc_{\ref{8.3}}(m) \geq 9 \cdot (m-1) \cdot 4^m$ and the proof is complete.
\end{proof}

\begin{Theorem}
\label{8.4}
There exists $f \in C_c(\Rm)$ with the following properties.
\begin{enumerate}
\item[(A)] There is a bounded Lipschitz vector field $v : \Rm \to \Rm$ whose distributional divergence coincides with $f$.
\item[(B)] There is no member of $C^1(\Rm;\Rm)$ whose distributional divergence coincides with $f$.
\end{enumerate}
\end{Theorem}

\begin{Remark}
\label{8.5}
If $w : \Rm \to \Rm$ is locally Lipschitz and $f : \Rm \to \R$ is locally Lebesgue-summable then $f$ is the distributional divergence of $w$ if and only if $(\rmdiv w)(x) = f(x)$ for $\ssfL^m$-a.e. $x \in \Rm$.
This yields equivalent formulations of conclusion (A) and (B) above.
\end{Remark}

\begin{proof}
{\bf (i)} 
For each $j \in \N_+$ choose $a_j \in \Rm$ with $|a_j|_2=2^{-j}$ and define $r_j = \frac{1}{16} \cdot 3^{-j}$.
Notice that $|a_j-a_k|_2 \geq \frac{1}{2} \cdot 2^{-\min\{j,k\}}$ whenever $j \neq k$ and that the balls $B(a_j,4\cdot r_j)$, $j \in \N_+$, are pairwise disjoint.
Choose also a decreasing sequence $\la \veps_j \ra_j$ of positive real numbers, converging to 0, and so that $\veps_1 \leq \frac{1}{8}$.
Applying corollary \ref{8.3} to $\veps_j$ provides us with $\hat{u}_j \in \calD_\#(\Rm)$ and $\hat{v}_j \in C^\infty(\Rm;\Rm)$ satisfying the conclusions $(A),\ldots,(F)$ there.
We define $u_j \in \calD_\#(\Rm)$ and $v_j \in C^\infty(\Rm;\Rm)$ by means of translation and rescaling as follows:
\begin{equation*}
u_j(x) = r_j^{-m} \cdot \hat{u}_j \left( \frac{x-a_j}{r_j}\right) \quad\text{and}\quad v_j(x) = r_j \cdot \hat{v}_j \left( \frac{x-a_j}{r_j}\right).
\end{equation*}
Abbreviate $f_j = \rmdiv v_j \in C^\infty(\Rm)$.
The conclusions of corollary \ref{8.3} become the following: 
\begin{enumerate}
\item $\rmspt u_j \subset U(a_j,r_j)$ and $\lno u_j \rno = \lno \hat{u}_j \rno \leq 2$ (by \ref{4.3}(B));
\item $\rmspt f_j \subset U(a_j,2 \cdot r_j)$;
\item $\|f_j\|_\infty \leq \veps_j$; 
\item $\displaystyle \int_{\Rm} u_j \cdot f_j \,d\ssfL^m \geq \frac{3}{8}$;
\item $\|v_j\|_\infty \leq \bc_{\ref{8.3}}(m) \cdot r_j$ and $\max_i \|\nabla (v_j \ip e_i)\|_\infty \leq \bc_{\ref{8.3}}(m)$;
\item For all $x \in \Rm \setminus B(a_j,4\cdot r_j)$ and all $i \in \{1,\ldots,m\}$ we have $\displaystyle |\nabla (v_j \ip e_i)(x)|_2 \leq \frac{\bc_{\ref{8.3}}(m) \cdot r_j^m}{|x-a_j|_2^m}$.
\end{enumerate}
\par 
{\bf (ii)}
We define $f : \Rm \to \R$ to coincide with $f_j$ in $U(a_j,2 \cdot r_j)$ for every $j \in \N_+$ and to vanish identically outside the union of these balls. 
Since these balls are pairwise disjoint, the definition makes sense and, in fact, $f = \sum_{j=1}^\infty f_j$, uniformly, since $\left\| f - \sum_{j=1}^n f_j \right\|_\infty \leq \veps_{n+1} \to 0$ as $n \to \infty$.
This proves at once that $f$ is continuous.
Also, notice that $\rmspt f \subset B(0,2)$.
Thus, $f \in C_c(\Rm)$.
\par 
{\bf (iii)}
A vector field $v : \Rm \to \Rm$ is well-defined by the formula $v = \sum_{j=1}^\infty v_j$.
In fact, the series converges normally, since $\sum_{j =1}^\infty \|v_j\|_\infty \leq \bc_{\ref{8.3}}(m) \sum_{j=1}^\infty r_j < \infty$.
In particular, $v$ is continuous and bounded.
We shall now establish that $v$ is Lipschitz.
To this end, we shall prove that the sequence $\la w_n \ra_n$ of partial sums $w_n = \sum_{j=1}^n v_j$ is equi-Lipschitz, equivalently, that each sequence $\la w_n \ip e_i \ra_n$, $i \in \{1,\ldots,m\}$, of coordinate functions of $w_n$ is equi-Lipschitz.
Fix $n \in \N_+$ and $x \in \Rm$.
Choose $k=k(x,n) \in \{1,\ldots,n\}$ such that $|x-a_k|_2=\min_{1\leq j\leq n}|x-a_j|_2$.
Recall that $|\nabla(v_{k} \ip e_i)(x)|_2 \leq \bc_{\ref{8.3}}(m)$. 
For $j \in \{1,\ldots,n\}\setminus\{k\}$ we have
\begin{equation*}
|x-a_j|_2 \geq \frac{1}{2} \cdot |a_k-a_j|_2 \geq \frac{1}{4} \cdot 2^{- \min \{k,j \}} \geq \frac{1}{4} \cdot 2^{-j} > 4 \cdot \frac{1}{16} \cdot 3^{-j} = 4 \cdot r_j.
\end{equation*} 
Thus, for such $j$, we have $x \not \in B(a_j , 4 \cdot r_j)$ and $\frac{1}{|x-a_j|_2} \leq 2^{j+2}$.
We infer that
\begin{equation*}
|\nabla(v_j \ip e_i)(x)|_2 \leq \frac{\bc_{\ref{8.3}}(m) \cdot r_j^m}{|x-a_j|_2^m} \leq \bc_{\ref{8.3}}(m) \cdot r_j^m \cdot \left( 2^{j+2} \right)^m = \bc_{\ref{8.3}}(m) \cdot 4^{-m} \cdot \left( \frac{2}{3}\right)^{j\cdot m}.
\end{equation*}
It follows that 
\begin{equation*}
|\nabla(w_n \ip e_i)(x)|_2 \leq \bc_{\ref{8.3}}(m) + \sum_{\substack{j=1 \\ j \neq k}}^n |\nabla (v_j \ip e_i)(x)|_2
\leq \bc_{\ref{8.3}}(m) \left( 1 + 4^{-m} \sum_{j=1}^\infty \left( \frac{2}{3} \right)^{jm} \right) =: \lambda.
\end{equation*}
Since $x \in \Rm$ is arbitrary, the gradient bound implies that $\rmLip (w_n \ip e_i) \leq \lambda$.
\par 
{\bf (iv)}
Recall that $\sum_{j=1}^n f_j \to f$ as $n \to \infty$ uniformly, hence, also in the sense of distributions.
As $\sum_{j=1}^n f_j = \sum_{j=1}^n \rmdiv v_j = \rmdiv \left( \sum_{j=1}^n v_j \right) = \rmdiv w_n$ and $w_n \to v$ distributionally, we conclude that $\sum_{j=1}^n f_j \to \rmdiv v$ as $n \to \infty$ distributionally.
Thus, $\rmdiv v = f$.
The proof of conclusion (A) is, therefore, complete.
\par 
{\bf (v)}
It remains to establish conclusion (B).
Arguing by contradiction, assume if possible that $f = \rmdiv w$ for some $w \in C^1(\Rm;\Rm)$.
Applying theorem \ref{3.3} with $\veps = \frac{1}{8}$ and $k=2$ would yield $\theta > 0$ such that for all $u \in \calD_\#(\Rm)$ if $\rmspt u \subset B(0,2)$ then 
\begin{equation}
\label{eq.8.4.contradiction}
\left| \int_{\Rm} u \cdot f \,d\ssfL^m \right| \leq \theta \cdot \bG(u) + \frac{1}{8} \cdot \lno u \rno.
\end{equation}
Upon recalling that $u_j \in \calD_\#(\Rm)$ and $\rmspt u_j \subset U(a_j,r_j) \subset B(0,2)$, we may apply \eqref{eq.8.4.contradiction} with $u=u_j$.
Noting that $\lno u_j \rno \leq 2$ and $\bG(u_j) = r_j \cdot \bG(\hat{u}_j) \leq r_j \cdot  \bc_{\ref{4.5}}(m,2) \cdot 2$, by \ref{4.3}(A) and theorem \ref{4.5}, we obtain
\begin{equation*}
\left| \int_{\Rm} u_j \cdot f \,d\ssfL^m \right| \leq \theta \cdot r_j \cdot \bc_{\ref{4.5}}(m,2) \cdot 2 + \frac{1}{4},
\end{equation*}
for all $j \in \N_+$.
As $f = f_j$ on $U(a_j,2 \cdot r_j)$ and $\rmspt u_j \subset U(a_j , 2 \cdot r_j)$, we have $\int_{\Rm} u_j \cdot f \,d\ssfL^m = \int_{\Rm} u_j \cdot f_j \,d\ssfL^m$.
Letting $j \to \infty$ yields a contradiction with {\bf (i)}(4).
\end{proof}

\section{$C(\Rm) \cap W^{1,m}_{\rmloc}(\Rm) \subset \rmdiv(C^1(\Rm;\Rm))$}

\begin{Empty}[The curl-div operator]
\label{9.1}
We start by recalling facts about cocancelling operators introduced by Van Schaftingen \cite{VSJ.13}.
Given two finite dimensional vector spaces $V$ and $W$, a linear differential operator homogeneous of degree $k \in \N_+$ with constant coefficients is a map $L(D) : \calD(\Rm;V) \to \calD(\Rm;W)$ such that 
\begin{equation*}
L(D)(u) = \sum_{\substack{\alpha \in \N^m \\ |\alpha|=k}} L_\alpha(D^\alpha u),
\end{equation*}
where as usual $D^\alpha = D_m^{\alpha_m} \circ \cdots \circ D_1^{\alpha_1}$ and $L_\alpha \in \rmHom(V,W)$.
Its symbol is defined to be
\begin{equation*}
\bbL : \Rm \to \rmHom(V,W) : \xi \mapsto \sum_{\substack{|\alpha|=k}} \xi^\alpha \cdot L_\alpha,
\end{equation*}
where as usual $\xi^\alpha = \xi_1^{\alpha_1} \cdots \xi_m^{\alpha_m}$.
We say that $L(D)$ is cocancelling if 
\begin{equation*}
\bigcap_{\xi \in \Rm \setminus \{0\}} \ker \bbL(\xi) = \{0\}.
\end{equation*}
Our interest in this notion is \cite[Theorem 1.4]{VSJ.13} that will be used in the proof of \ref{9.3}(iv).
\par 
We define the curl-div operator $L(D)$ as above by setting $k=2$, $V = W = M_m(\R)$, and for $A \in \calD(\Rm;M_m(\R))$ the matrix $L(D)(A) \in M_m(\R)$ is defined by its entries
\begin{equation*}
\big( L(D)(A) \big)_{ij} = D_j( \rmdiv A_i ) - D_i( \rmdiv A_j),
\end{equation*}
where $A_i$ and $A_j$ are, respectively, the $i^{th}$ and $j^{th}$ rows of $A$.
\par 
Decompose $L(D) = L_1(D) - L_2(D)$ in the obvious way.
Thus, the $i^{th}$ row of $L_1(D)(A)$ is the gradient of the divergence of the $i^{th}$ row of $A$.
In order to determine the symbol of $L_1(D)$ we introduce the linear maps $C_{lk} : M_m(\R) \to M_m(\R)$ so that $C_{lk}(A)$ is the matrix whose $l^{th}$ column is the $k^{th}$ column of $A$ and all other entries vanish, \ie $C_{lk}(A) = (A \cdot e_k) \otimes e_l$.
In other words, $\big( C_{lk}(A)\big)_{ij} = \delta_{jl} \cdot A_{ik}$.
Observe that
\begin{equation*}
\big( L_1(D)(A) \big)_{ij} = \sum_{k=1}^m (D_j \circ D_k)(A_{ik}) = \sum_{k,l=1}^m \delta_{jl} \cdot (D_l \circ D_k)(A_{ik}) = \sum_{k,l=1}^m (D_l \circ D_k)\big( C_{lk}(A) \big)_{i,j} .
\end{equation*}
Accordingly,
\begin{equation*}
L_1(D)(A) = \sum_{k,l=1}^m (D_l \circ D_k)\big( C_{lk}(A) \big).
\end{equation*}
It follows that
\begin{equation*}
\bbL_1(\xi)(A) = \sum_{k,l=1}^m \xi_l \cdot \xi_k \cdot (A \cdot e_k) \otimes e_l = (A \cdot \xi) \otimes \xi.
\end{equation*}
Similarly, $\bbL_2(\xi)(A) = \xi \otimes (A \cdot \xi)$.
We conclude that 
\begin{equation}
\label{eq.9.1}
\bbL : \Rm \to \rmHom(M_m(\R),M_m(\R)) : \xi \mapsto \big[ A \mapsto (A \cdot \xi) \otimes \xi - \xi \otimes (A \cdot \xi) \big].
\end{equation}
\begin{enumerate}
\item[(A)] {\it $\displaystyle \cap_{\xi \in \Rm \setminus \{0\}} \ker \bbL(\xi) = \R \cdot I_m.$}
\end{enumerate}
\par 
{\it Proof.}
Assume $A \in M_m(\R)$ belongs to $\ker \bbL(\xi)$ for all $\xi \in \Rm \setminus \{0\}$.
Observe that $\bbL(\xi)(A) = 0$ iff $(A \cdot \xi) \otimes \xi = \xi \otimes (A \cdot \xi)$ iff $\xi$ is an eigenvector of $A$.
This holds for all $\xi \neq 0$ iff $A$ is a scalar multiple of the identity matrix $I_m$.\cqfd
\par 
Recalling the notation $M_{m,0}(\R)$ of trace-free matrices from \ref{np.matrices}, we consider the curl-div restricted operator
\begin{equation*}
\widetilde{L(D)} : \calD(\Rm;M_{m,0}(\R)) \to \calD(\Rm;M_m(\R)) : A \mapsto L(D)(A).
\end{equation*}
\begin{enumerate}
\item[(B)] {\it $\widetilde{L(D)}$ is cocancelling.}
\end{enumerate}
\par 
{\it Proof.}
This is an immediate consequence of (A) and $\ker(\rmtrace) \cap (\R \cdot I_m) = \{0\}$.\cqfd
\par 
We shall also need the following two elementary observations.
\begin{enumerate}
\item[(C)] {\it Assume that $v_i \in \calD(\Rm;\Rm)$, $i \in \{1,\ldots,m\}$, and there exists $u \in \calD_\#(\Rm)$ such that $\rmdiv v_i = D_iu$ for all $i$. If the rows of $A \in \calD(\Rm;M_m(\R))$ are $v_1,\ldots,v_m$ then $L(D)(A)=0$.}
\end{enumerate}
\par 
{\it Proof.}
This follows at once from the Clairaut-Schwarz-Young theorem.\cqfd
\begin{enumerate}
\item[(D)] {\it $L(D)(\vphi \cdot I_m) = 0$ for all $\vphi \in \calD(\Rm)$.}
\end{enumerate}
\par 
{\it Proof.}
Again, $D_jD_i \vphi = D_iD_j\vphi$.\cqfd
\end{Empty}

\begin{Proposition}
\label{9.2}
If $u \in \calD_\#(\Rm)$, $k \in \N_+$, and $\rmspt u \subset B(0,k)$ then there exists $v \in \calD(\Rm;\Rm)$ such that $\rmdiv v = u$, $\|v\|_{L_1} \leq 3 \cdot \bG(u)$, and $\rmspt v \subset B(0, k+1)$.
\end{Proposition}

\begin{proof}
{\bf (i)}
We may assume without loss of generality that $u \neq 0$, hence, also that $\bG(u) > 0$.
Put $\veps = \bG(u)$ and choose $0 < \delta \leq 1$ small enough for $\delta \cdot \|u\|_{L_1} < \bG(u)$.
We infer from \ref{4.2}(C) that there exists $S \in \calD_1(\Rm)$ such that $\partial S = \iota(u)$ and $\bM(S) < \veps + \bG(\iota(u)) = 2 \cdot \bG(u)$.
Letting $\pi : \Rm \to \Rm$ denote the nearest point projection onto $B(0,k)$, we infer from \cite[4.1.15 end on p. 372]{GMT} that $\partial \pi_\#S = \pi_\# \partial S = \pi_\# \iota(u) = \iota(u)$, since $\rmspt \iota(u) \subset B(0,k)$. 
As $\bM(\pi_\#S) \leq \bM(S)$, replacing $S$ by $\pi_\#S$ if necessary, we may assume that $\rmspt S \subset B(0,k)$.
There exists $w_\delta \in \calD(\Rm;\Rm)$ such that $S * \Phi_\delta = -\ssfL^m \wedge w_\delta$ and $\rmspt w_\delta \subset B(0,k+1)$ and one routinely checks that $\rmdiv w_\delta = u * \Phi_\delta$.
Furthermore, $\|w_\delta\|_{L_1} = \bM(S * \Phi_\delta) \leq \bM(S) \leq 2 \cdot \bG(u)$.
\par 
{\bf (ii)}
Next, we define $w \in \calD(\Rm;\Rm)$ by the formula
\begin{equation*}
w(x) = \int_{\Rm \times [0,1]} \Phi_\delta(y) \cdot u(x-t \cdot y) \cdot y \, d\ssfL^{m+1}(y,t).
\end{equation*}
Observe that $\rmspt w \subset B(0,k+1)$, that
\begin{equation*}
\begin{split}
(\rmdiv w)(x) & = \int_{\Rm \times [0,1]} \Phi_\delta(y) \cdot \big( (\nabla u)(x-t \cdot y) \ip y \big) \, d\ssfL^{m+1}(y,t) \\
& = -\int_{\Rm} \Phi_\delta(y) \left( \int_0^1 \frac{d}{dt} u(x-t \cdot y) \,d\ssfL^1(t)\right) \,d\ssfL^m(y) \\
& = \int_{\Rm} \Phi_\delta(y) \cdot \big( u(x) - u(x-y) \big)\,d\ssfL^m(y) \\
& = u(x) - (u * \Phi_\delta)(x),
\end{split}
\end{equation*}
and that $\|w\|_{L_1} \leq \delta \cdot \|u\|_{L_1} < \bG(u)$.
\par 
{\bf (iii)}
Finally, we define $v = w_\delta + w \in \calD(\Rm;\Rm)$, so that $\rmdiv v = \rmdiv w_\delta + \rmdiv w = u$, $\|v\|_{L_1} \leq \|w_\delta\|_{L_1} + \|w\|_{L_1} \leq 3 \cdot \bG(u)$, and $\rmspt v \subset B(0,k+1)$.
This completes the proof of the proposition.
\end{proof}

\begin{Theorem}
\label{9.3}
There exist $\bc_{\theTheorem}(m,k) > 0$ with the following property.
Let $f \in C(\Rm)$ and assume that $D_i f \in L_{m,\rmloc}(\Rm)$ for all $i \in \{1,\ldots,m\}$.
For every $u \in \calD_\#(\Rm)$ and $k \in \N_+$ such that $\rmspt u \subset B(0,k)$ we have
\begin{equation*}
\left| \int_{\Rm} u \cdot f \,d\ssfL^m \right| \leq \bc_{\theTheorem}(m,k) \cdot \left( \|f\|_{\infty,B(0,k+2)} + \|\nabla f\|_{L_m,B(0,k+2)} \right) \cdot \lno u \rno.
\end{equation*}
\end{Theorem}

\begin{proof}
{\bf (i)}
Applying proposition \ref{9.2} to each $D_iu$ (in place of $u$), $i \in \{1,\ldots,m\}$, we obtain $v_i \in \calD(\Rm;\Rm)$ such that $\rmdiv v_i = D_iu$, $\|v_i\|_{L_1} \leq 3 \cdot \bG(D_iu)$, and $\rmspt v_i \subset B(0,k+1)$.
Note that $\sum_{i=1}^m \|v_i\|_{L_1} \leq 3 \cdot \lno u \rno$.
\par 
{\bf (ii)}
We shall prove the theorem first in case $f \in \calD(\Rm)$.
Choose $\psi \in \calD(\Rm)$ such that $\ind_{B(0,k+1)} \leq \psi \leq \ind_{B(0,k+2)}$ and $\|\nabla\psi\|_\infty \leq 2$. Define $g = \Gamma * (\psi \cdot f) \in C^\infty(\Rm)$, so that $\triangle g = \psi \cdot f = f$ on a neighbourhood of $B(0,k+1)$.
Integrating by parts, we infer that
\begin{multline*}
\int_{\Rm} u \cdot f \, d\ssfL^m = \int_{\Rm} u \cdot \triangle g \,d\ssfL^m = - \int_{\Rm} (\nabla u) \ip (\nabla g) \,d\ssfL^m = - \sum_{i=1}^m \int_{\Rm} D_iu \cdot D_i g \,d\ssfL^m \\
= - \sum_{i=1}^m \int_{\Rm} (\rmdiv v_i) \cdot D_ig \, d\ssfL^m = \sum_{i=1}^m \int_{\Rm} v_i \ip (\nabla D_i g) \,d\ssfL^m = \int_{\Rm} A : D^2g \,d\ssfL^m,
\end{multline*}
where $A \in \calD(\Rm;M_m(\R))$ is the matrix field whose $i^{th}$ row is $v_i$.
Defining $A^\circ := A - \frac{1}{m} \cdot (\rmtrace A) \cdot I_m$ we note that $A^\circ \in \calD(\Rm;M_{m,0}(\R))$ and $A = \frac{1}{m} \cdot (\rmtrace A) \cdot I_m + A^\circ$.
Therefore,
\begin{equation*}
\int_{\Rm} u \cdot f \,d\ssfL^m = \frac{1}{m} \cdot \int_{\Rm} (\rmtrace A) \cdot I_m : D^2g \,d\ssfL^m + \int_{\Rm} A^\circ : D^2 g \,d\ssfL^m = \rmI + \rmII,
\end{equation*}
where $\rmI$ and $\rmII$ are defined by the above equation.
\par 
{\bf (iii)}
In order to estimate $|\rmI|$ from above we note that $I_m : D^2g = \rmtrace(D^2g) = \triangle g = f$ on the support of $A$, thus,
\begin{multline*}
|\rmI| = \frac{1}{m} \cdot \left| \int_{\Rm} (\rmtrace A) \cdot f \,d\ssfL^m \right| \leq
\|\rmtrace A\|_{L_1} \cdot \|f\|_{\infty,B(0,k+1)}  \\ \leq 
\left( \sum_{i=1}^m \|v_i\|_{L_1} \right) \cdot \|f\|_{\infty,B(0,k+1)}
\leq 3 \cdot \lno u \rno \cdot \|f\|_{\infty,B(0,k+1)}.
\end{multline*}
\par 
{\bf (iv)}
In order to obtain an upper bound for $|\rmII|$ we shall apply \cite[Theorem 1.4]{VSJ.13} to the operator $\widetilde{L(D)}$ (recall \ref{9.1}) and the trace-free matrix field $A^\circ$.
To this end, we observe that $L(D)(A) = 0$, by \ref{9.1}(C), and that $L(D)(\rmtrace(A) \cdot I_m) = 0$, by \ref{9.1}(D), whence, $\widetilde{L(D)}(A^\circ) = L(D)(A) - \frac{1}{m} \cdot L(D)((\rmtrace A) \cdot I_m) = 0$.
Accordingly,
\begin{equation*}
|\rmII| = \left| \int_{\Rm} A^\circ : D^2g \,d\ssfL^m \right| = \left| \int_{\Rm} A^\circ : (\psi \cdot D^2g) \,d\ssfL^m \right| \leq \bc_{\mathrm{VS}}(m) \cdot \|A^\circ\|_{L_1} \cdot \|D(\psi \cdot D^2g)\|_{L_m},
\end{equation*} 
where $\bc_{\mathrm{VS}}(m)$ is the dimensional constant in \cite[Theorem 1.4]{VSJ.13}.
Recall that $\|A^\circ\|_{L_1} \leq 2 \cdot \|A\|_{L_1} = 2 \cdot \sum_{i=1}^m \|v_i\|_{L_1} \leq 6 \cdot \lno u \rno$, according to \ref{np.matrices}.
Regarding the last factor in the inequality above, the global Calder\'on-Zygmund inequality (see \cite[Chapter II, \S 4, 4.2(b)]{STEIN.70}) applied on $\Rm$ to $\psi \cdot f$ and its first derivatives gives
\begin{equation*}
\begin{split}
\|D(\psi \cdot D^2g)\|_{L_m(\Rm)}
& \leq \|D^2(\nabla g)\|_{L_m(\Rm)} + 2 \cdot \|D^2g\|_{L_m(\Rm)} \\
& \leq \bc_{\mathrm{CZ}}(m) \cdot \big(\|\nabla(\psi \cdot f)\|_{L_m(\Rm)} + \|\psi \cdot f\|_{L_m(\Rm)}\big) \\
& \leq \bc'_{\mathrm{CZ}}(m,k) \cdot \big(\|\nabla f\|_{L_m(B(0,k+2))} + \|f\|_{\infty,B(0,k+2)}\big),
\end{split}
\end{equation*}
where we used $\triangle g = \psi \cdot f$, $\|\psi\|_\infty \leq 1$, $\|\nabla\psi\|_\infty \leq 2$, and $\rmspt\psi \subset B(0,k+2)$.
We conclude that
\begin{equation*}
|\rmII| \leq \bc(m,k) \cdot \left( \|f\|_{\infty,B(0,k+2)} + \|\nabla f\|_{L_m(B(0,k+2))}\right) \cdot \lno u \rno
\end{equation*}
for some appropriate $\bc(m,k) > 0$.
\par 
{\bf (v)}
Combining the estimates in {\bf (ii)}, {\bf (iii)}, and {\bf (iv)} yields the conclusion of the theorem in case $f \in \calD(\Rm)$.
In the general case, we choose $\psi \in \calD(\Rm)$ which equals 1 in a neighbourhood of $B(0,k+2)$ and we apply the particular case to $\Phi_\delta * (\psi \cdot f)$ for each $\delta > 0$.
It remains to notice that $\int_{\Rm} u \cdot f \, d\ssfL^m = \lim_{\delta \to 0^+} \int_{\Rm} u \cdot [\Phi_\delta * (\psi \cdot f)] \, d\ssfL^m$, $\|f\|_{\infty,B(0,k+2)} = \lim_{\delta \to 0^+} \|\Phi_\delta *( \psi \cdot f) \|_{\infty,B(0,k+2)}$, and $\|\nabla f\|_{L_m(B(0,k+2))} = \lim_{\delta \to 0^+} \|\nabla [ \Phi_\delta * (\psi \cdot f)]\|_{L_m(B(0,k+2))}$.
This completes the proof.
\end{proof}

\begin{Theorem}
\label{9.4}
Let $f \in C(\Rm)$ and assume that $D_if \in L_{m,\rmloc}(\Rm)$ for all $i \in \{1,\ldots,m\}$.
Then there exists $v \in C^1(\Rm;\Rm)$ such that $\rmdiv v = f$.
In fact, for each $k \in \N_+$ such that $\int_{\Rm} u \cdot f \,d\ssfL^m \neq 0$ for some $u \in \calD_\#(\Rm)$ with $\rmspt u \subset B(0,k)$ there exists $v \in C^1(\Rm;\Rm)$ (possibly depending on $k$) such that $\rmdiv v = f$ and
\begin{equation*}
\bq_k(v) \leq \left( 4 \cdot \bc_{\ref{9.3}}(m,k+1) +1 \right)\cdot \left( \|f\|_{\infty,B(0,k+3)} + \|\nabla f\|_{L_m,B(0,k+3)} \right).
\end{equation*}
\end{Theorem}

\begin{proof}
{\bf (i)}
According to theorem \ref{6.8}, it suffices first to establish that
\begin{equation*}
(\forall \veps > 0)(\forall k \in \N_+)(\exists \theta > 0)(\forall u \in \calD_\#(\Rm)): \rmspt u \subset B(0,k) \Rightarrow \left| \int_{\Rm} u \cdot f \, d\ssfL^m \right| \leq \theta \cdot \bG(u) + \veps \cdot \lno u \rno.
\end{equation*}
Fix $\veps > 0$ and $k \in \N_+$.
Let $\psi \in \calD(\Rm)$ such that $B(0,k) \subset \{ \psi = 1 \}$.
For every $\delta > 0$, applying theorem \ref{9.3} to $f - \Phi_\delta * f$ we obtain, for each $u \in \calD_\#(\Rm)$ supported in $B(0,k)$,
\begin{equation*}
\begin{split}
\left| \int_{\Rm} u \cdot f \, d\ssfL^m \right| & \leq \left| \int_{\Rm} u \cdot \psi \cdot (f * \Phi_\delta) \, d\ssfL^m \right| + \left| \int_{\Rm} u \cdot (f - f *\Phi_\delta) \,d\ssfL^m \right| \\
& \leq \|\nabla[ \psi \cdot (f * \Phi_\delta)]\|_\infty \cdot \bG(u) \\
& \qquad\quad + \bc_{\ref{9.3}}(m,k) \cdot \left( \|f - f * \Phi_\delta\|_{\infty,B(0,k+2)} + \|\nabla f - (\nabla f) * \Phi_\delta\|_{L_m,B(0,k+2)}\right) \cdot \lno u \rno.
\end{split}
\end{equation*}
Choosing $\delta > 0$ so small that
\begin{equation*}
\bc_{\ref{9.3}}(m,k) \cdot \left( \|f - f * \Phi_\delta\|_{\infty,B(0,k+2)} + \|\nabla f - (\nabla f) * \Phi_\delta\|_{L_m,B(0,k+2)}\right) < \veps
\end{equation*}
and letting
\begin{equation*}
\theta = 1 + \|\nabla[ \psi \cdot (f * \Phi_\delta)]\|_\infty
\end{equation*}
proves the requested condition. By theorem \ref{6.8}, there exists $w \in C^1(\Rm;\Rm)$ such that
\begin{equation*}
\int_{\Rm} u \cdot (\rmdiv w-f) \,d\ssfL^m=0 \qquad\text{for all }u \in \calD_\#(\Rm).
\end{equation*}
By the constancy theorem, $\rmdiv w-f=c$ for some $c \in \R$. Thus $v=w-\frac{c}{m}\cdot\rmid_{\Rm}$ satisfies $\rmdiv v=f$, proving the first conclusion.
\par
{\bf (ii)}
We now turn to proving the estimate. 
For $k \in \N_+$, abbreviate
\begin{equation*}
M_k=\|f\|_{\infty,B(0,k+3)}+\|\nabla f\|_{L_m,B(0,k+3)}.
\end{equation*}
Let $F : \calD_\#(\Rm) \to \R$ be defined by $F(u)=\int_{\Rm}u\cdot f\,d\ssfL^m$ and let $\hat F \in X[\calG_{\lno\cdot\rno}]^*$ be its unique continuous extension, according to theorem \ref{6.8}.
Recall from \ref{6.5} that, for $T \in X_k$, there are $u_j \in \calD_\#(\Rm)$ supported in $B(0,k+1)$ such that $\bG(\iota(u_j)-T)\to0$ and $\sup_j\lno u_j\rno\leq\lno T\rno$. Applying theorem \ref{9.3} to $u_j$ and letting $j\to\infty$, we obtain
\begin{equation*}
|\la T,\hat F\ra|\leq\bc_{\ref{9.3}}(m,k+1)\cdot M_k\cdot\lno T\rno.
\end{equation*}
Consequently, writing
\begin{equation*}
\vvvert \hat{F} \vvvert_k=\sup\big\{|\la T,\hat F\ra| \mid T\in X_k\text{ and }\lno T\rno\leq1\big\},
\end{equation*}
we have $0 < \vvvert \hat{F} \vvvert_k \leq\bc_{\ref{9.3}}(m,k+1)\cdot M_k$, where the first inequality is in the assumption of the theorem.
It follows from \cite[Theorem 8.1(C3)]{DEP.26c}, applied with $\veps=1$, that there exists $w \in \bE^1(\Rm)$ such that $\bD(w) = \hat{F}$ and $\bq_k(w) \leq 2 \cdot \vvvert \hat{F} \vvvert_k \leq 2 \cdot \bc_{\ref{9.3}}(m,k+1) \cdot M_k$.
As in {\bf (i)}, $\rmdiv w - f = c$ for some $c \in \R$ and $\rmdiv v = f$, where $v = w - \frac{c}{m} \cdot \rmid_{\Rm}$.
Notice that $|c| = |(\rmdiv w)(0) - f(0)| \leq m \cdot \bq_k(w) + M_k$.
Thus, $\bq_k(v) \leq \bq_k(w) + \frac{|c|}{m} \leq \left( 4 \cdot \bc_{\ref{9.3}}(m,k+1) + \frac{1}{m} \right) \cdot M_k$.
This completes the proof.
\end{proof}


\section*{Tool and computational resource disclosure}

No computer algebra system, proof assistant, or numerical computation was used to establish or verify the mathematical results of this paper.
Several large language models (LLMs) were used for editorial assistance and the critical examination of mathematical arguments.
In particular, an LLM proposed a proof of Theorem \ref{8.4} which was found to be incorrect and was discarded.
On the other hand, an LLM suggested the use of Van Schaftingen's theorem on cocancelling operators in the proof of Theorem \ref{9.4}.
The author takes full responsibility for the correctness of the results and their proofs.

\printbibliography




\end{document}